\documentclass[a4paper,12pt]{article}

\usepackage{ae,aeguill}
\usepackage{amsmath}
\usepackage{amsfonts}
\usepackage{amssymb}
\usepackage[T1]{fontenc}
\usepackage{geometry}
\usepackage{graphicx}
\usepackage{hyperref}
\usepackage{icomma}
\usepackage[latin1]{inputenc}
\usepackage{latexsym}
\usepackage[square]{natbib}
\usepackage{textcomp}
\usepackage[all]{xy}

\newcommand*\FlechBH[3]{\smash{\mathop{#1}\limits_{#2}^{#3}}}

\DeclareMathSymbol{\minus}{\mathord}{operators}{"2D}

\newtheorem{theoreme}{Theorem}
\newtheorem{lemme}[theoreme]{Lemma}
\newtheorem{proposition}[theoreme]{Proposition}
\newtheorem{remarque}{Remark}
\newtheorem{definition}{Definition}
\newtheorem{hypothese}{Hypothesis}

\newtheorem{extthm}{Theorem}

\title{Approximation of quasi-stationary distributions for 1-dimensional killed diffusions with unbounded drifts}

\author{Denis Villemonais}

\begin{document}

\maketitle

\begin{abstract}
 The long time behavior of an absorbed Markov process is well described by the limiting distribution of the process conditioned to not be killed when it is observed. Our aim is to give an approximation's method of this limit, when the process is a $1$-dimensional It\^o diffusion whose drift is allowed to explode at the boundary. In a first step, we show how to restrict the study to the case of a diffusion with values in a bounded interval and whose drift is bounded. In a second step, we show an approximation method of the limiting conditional distribution of such diffusions, based on a Fleming-Viot type interacting particle system. We end the paper with two numerical applications : to the logistic Feller diffusion and to the Wright-Fisher diffusion with values in $]0,1[$ conditioned to be killed at $0$.
\end{abstract}


\noindent\textit{MSC 2000 subject : }Primary 65C50, 60K35; secondary 60J60

\section{Introduction}

Let $(X_t)$ be a killed Markov process with law $\mathbb{P}$, taking its values in $E \cup \{\partial\}$, where $\partial$ is a cemetery point. We denote by $\tau_{\partial}=\inf\{t\geq0, X_t=\partial\}$ the killing time of $(X_t)$. A probability measure $\nu$ on $E$ is called a \textbf{quasi-stationary distribution} (QSD) if, for all $t\geq0$, the distribution of the process $X$, initially distributed with respect to $\nu$ and conditioned to be not killed before time $t$, is still $\nu$ at time $t$, that is
 $\mathbb{P}_{\nu}\left(X_t\in A|\tau_{\partial}> t\right)=\nu(A)$ for every $A\subset E$ and $t\geq 0$. Without loss of generality, we suppose that  $\partial$ is an absorbing point, so that $\{\tau_{\partial}>t\}=\{X_t\neq\partial\}$.

Let $\mu$ be a probability measure on $E$. If it exists and provided it is a probability, the limiting conditional distribution
\begin{equation*}
 \lim_{t\rightarrow+\infty}{\mathbb{P}_x(X_t\in . |X_t\neq\partial)}
\end{equation*}
is called the Yaglom limit for $\mu$, from the Russian Mathematician A.M. Yaglom. He showed in \cite{ya1} that the limiting conditional distribution of the number of descendants in the $n^{th}$ generation of a Galton-Watson process always exists in the subcritical case.

The existence or uniqueness of such invariant conditional distributions have been proved in a host of contexts. When $E$ is finite, it is proved in \cite{da1} that there exists a unique QSD $\nu$ and that the Yaglom limit converges to $\nu$ independently of the initial distribution. In \cite{ca2}, the case of a birth and death process on $\mathbb{N}$ is studied. For this process, the set of QSDs is either empty, or a singleton, or a continuum indexed by a real parameter and given by an explicit recursive formula. This is an exception : most of the known results on QSDs are related with existence or uniqueness problems. In \cite{fe2}, the existence of a quasi-stationary distribution for a continuous time Markov chain on $\mathbb{N}$ killed at $0$ is proved under conditions on moments of the killing time, using an original renewal dynamical approach. In \cite{co1}, the case of $1$-dimensional diffusion on $[0,+\infty[$ with $C^1$ drift and killed at $0$ are studied, with the assumption that $+\infty$ is a natural boundary. The dependence between the initial measure and the Yaglom limit is explored in \cite{ma1} (for a Brownian motion with constant drift killed at $0$) and \cite{ll1} (for the Orstein-Uhlenbeck process killed at $0$). In \cite{st2}, the case of $1$-dimensional diffusions with general killing on the interior of a given interval is investigated. In \cite{ca1}, the authors study the existence and uniqueness of the QSD for $1$-dimensional diffusions killed at $0$ and whose drift is allowed to explode at the boundary, which is the case under study in the present paper. See \cite{po1} for a regularly updated extensive bibliography on QSD.

In this paper we are concerned with $1$-dimensional It\^o diffusions with values in $]0,+\infty[\cup\{\partial\}$ killed at $0$ and defined by the stochastic differential equation
\begin{equation*}
\label{eq73}
dX_t=dB_t-q(X_t)dt,\ X_0=x>0,
\end{equation*}
 where $B$ is a standard $1$-dimensional Brownian motion and $q\in C^1(]0,+\infty[)$. In \cite{ca1}, the Yaglom limit of this process is studied and the authors give some conditions on the drift $q$, which are sufficient for the existence and the uniqueness of the QSD. In particular, they allow the drift to explode at the origin. As explained in the paper, this diffusion is closely related with some Markov mortality models. Such applications need the computation of the process QSD, but the tools used in \cite{ca1} are based on spectral theory's arguments and don't allow us to get explicit values. Our aim is to give an easily simulable approximation's method of this QSD.

The problem of QSD's approximation has been already explored in \cite{bu2}, \cite{gr1} when $E$ is a bounded open set of $\mathbb{R}^d$ and $X$ is a Brownian motion killed at the boundary of $E$. The authors proved an approximation's method exposed in \cite{bu1}, which is based on a Fleming-Viot type system of interacting particles whose number is going to infinity. In \cite{fe1}, it is proved that this method works well for a continuous time Markov chain in a countable state space under suitable assumption on the transition's rates (moreover, the existence of a QSD is a consequence of the approximation's method). New difficulties arise from our case with unbounded drift. For instance, the interacting particle process introduced in \cite{bu1} isn't necessarily well defined. To avoid this difficulty, we begin by proving that one can approximate our QSD by the QSD's of diffusions with bounded drifts.

\begin{figure}[htbp]
\begin{center}
\input{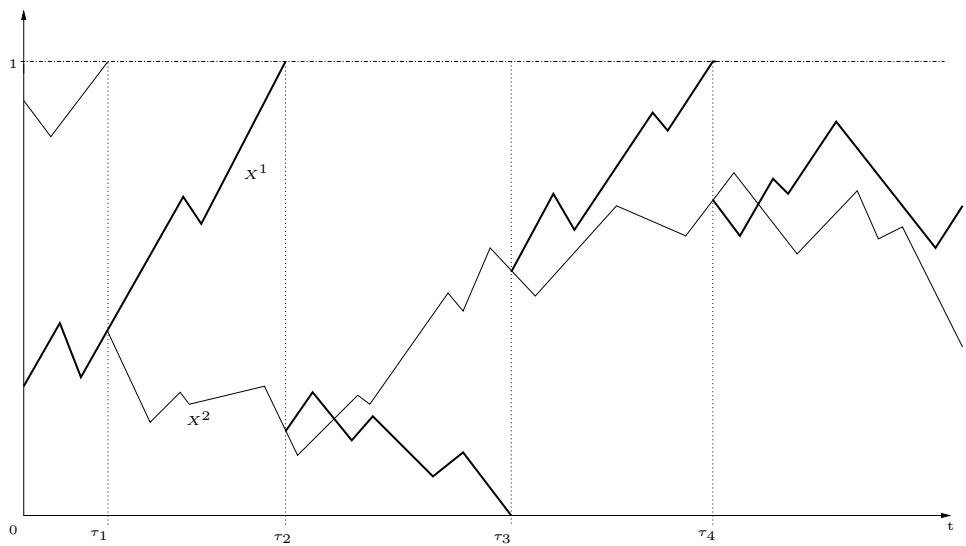}
\caption{The interacting particle system $(X^1,X^2)$}
\label{fig1}
\end{center}
\end{figure}

Let us denote by $\mathbb{P}^{\epsilon}$ the law of a diffusion with values in $]\epsilon,1/\epsilon[$, defined by the stochastic differential equation $dX_t=dB_t-q(X_t)dt$ and killed when it hits $\epsilon$ or $1/\epsilon$.
In \cite{pi1}, it is proved that the Yaglom limit associated with $\mathbb{P}^{\epsilon}$ exists and is its unique QSD. We will denote it by $\nu_{\epsilon}$.
In the first part of this paper, we give some conditions on $q\in C^1(]0,+\infty[)$ for the family $(\nu_{\epsilon})_{0<\epsilon\leq1/2}$ to be tight and to converge, when $\epsilon\rightarrow0$, to a QSD for the law $\mathbb{P}^0$.
We point out the fact that this result remains valid in the case of an unbounded drift diffusion with values in a bounded interval.
In a second part, we prove an approximation method for each probability measure $\nu_{\epsilon}$, based on the interacting process introduced in \cite{bu1}.
Fix $\epsilon>0$ and let us describe the interacting particle process of size $N\geq2$: each particle moves independently in $]\epsilon,1/\epsilon[$, each one with  law $\mathbb{P}^{\epsilon}$ until one of them hits the boundary. At this time, the killed particle jumps on the position of an other particle, chosen uniformly between the $N-1$ remaining one. Then the particles evolve independently, until one of them is killed and so on (see Figure \ref{fig1}).
One has to prove that the particles don't degenerate at the boundary. In \cite{bu2}, the authors prove a non-degeneracy result with arguments based on a construction of the $d$-dimensional Brownian motion due to It\^o, where $d\geq2$. It seems that this tool can't be easily generalized to other diffusions. To prove such results under our settings, we build an original coupling between the interacting particle process and an independent particle system of the same size. This coupling is valid for all drifted Brownian motions with continuous bounded drift, killed at the boundary of a bounded interval of $\mathbb{R}_+$. It will be used in each step of the proof.

We conclude the paper by two numerical applications. At first, we treat the case of the logistic Feller diffusion introduced in \cite{la1} and studied in \cite{ca1} with values in $]0,+\infty[$, driven by the stochastic differential equation
\begin{equation*}
 \label{eq74}
 dZ_t=\sqrt{Z_t}dB_t+(r Z_t-c Z_t^2)dt,\ Z_0=z>0,
\end{equation*}
where $B$ is a $1$-dimensional Brownian motion and $r, c$ are two positive constants. Clearly, $0$ is an absorbing state for this diffusion.
 In a second time we study in detail the case of the Wright-Fisher diffusion on $]0,1[$ conditioned to be killed at $1$ (see \cite{hu1}). This diffusion takes values in $]0,1[$, and is defined by
\begin{equation*}
 \label{eq75}
 dZ_t=\sqrt{Z_t(1-Z_t)}dB_t+(1-Z_t)dt,\ Z_0=z\in]0,1[,
\end{equation*}
where $B$ is a $1$-dimensional Brownian motion. This diffusion is absorbed at $1$.

\section{From unbounded drift to bounded drift}
\label{par1}
Let $\mathbb{P}^0_x$ be the law of a diffusion process taking its values in $]0,+\infty[\cup\{\partial\}$, killed when it hits $0$ and defined by the stochastic differential equation (SDE)
\begin{equation*}
dX_t=dB_t-q(X_t)dt,\ X_0=x>0,
\end{equation*}
where $B$ is a $1$-dimensional Brownian motion. The drift $q$ is taken in the set of real valued continuously differentiable functions $C^1(]0,+\infty[)$. We denote by $L^0$ the infinitesimal generator associated with $\mathbb{P}^0$.

We define, $\forall x\in]0,+\infty[$,
\begin{equation*}
Q(x)=\int_1^{x}{q(y)dy},
\end{equation*}
\begin{equation*}
 d\mu(x)=e^{-2Q(x)}dx,
\end{equation*}
and
\begin{equation*}
 W(x)=q(x)^2-q'(x).
\end{equation*}

For all $\epsilon\in]0,1/2[$, we define $\mathbb{P}^{\epsilon}_x$ as the law of the diffusion taking its values in $]\epsilon,1/\epsilon[$, defined by the SDE 
\begin{equation*}
 dX_t=dB_t-q(X_t)dt,\ X_0=x\in]\epsilon,1/\epsilon[
\end{equation*}
and killed when it hits the boundary $\{\epsilon,1/\epsilon\}$. Let $L^{\epsilon}$ be the infinitesimal generator $1/2\Delta-q\nabla$ with the Dirichlet boundary condition on $\{\epsilon,1/\epsilon\}$. $-L^{\epsilon}$ has a simple real eigenvalue $\lambda_{\epsilon}$ (see \cite[Theorem KR]{pi1}) at the bottom of its spectrum. The corresponding eigenfunction $\eta_{\epsilon}$ is positive and belongs to $C^2([\epsilon,1/\epsilon])$. We choose it so that
\begin{equation}
 \label{eq70}
 \int_{\epsilon}^{1/\epsilon}{\eta_{\epsilon}(x)^2 d\mu(x)}=1.
\end{equation}
Let us recall some results of \cite{pi1}:
\begin{extthm}[Pinsky (1985)]
 \label{th9}
 The Yaglom limit associated with $\mathbb{P}^{\epsilon}$ exists for all initial distributions $\delta_x$, $x\in]\epsilon,1/\epsilon[$, and doesn't depend on $x$. This limit is a QSD, which we denote by $\nu_{\epsilon}$. Furthermore, we have
\begin{equation}
\label{eqb1}
 d\nu_{\epsilon}(x)=\frac{\eta_{\epsilon}(x)d\mu(x)}{\int_{\epsilon}^{1/\epsilon}{\eta_{\epsilon}(x)d\mu(x)}}.
\end{equation}
\end{extthm}
In fact, $\nu_{\epsilon}$ is the unique QSD of the process, as proved in Lemma \ref{le3} below, but we won't use it in this section. The aim of this section is to study the asymptotic behaviour of $(\nu_{\epsilon})$ when $\epsilon$ goes to $0$.

From now, ${\cal M}_1(]0,+\infty[)$ denotes the space of probability measures on $]0,+\infty[$ equipped with the weak topology. The following hypotheses arise naturally in the proof of Theorem \ref{th1}, which is based on a compactness-uniqueness method. 
\begin{hypothese}[H1]
 $W$ is bounded below by $-C$, where $C$ is a positive constant. Moreover, $W(x)\rightarrow+\infty$ when $x\rightarrow\infty$.
\end{hypothese}
\begin{hypothese}[H2]
 \begin{equation*}
  \int_1^{+\infty}{e^{-2Q(x)}dx}<+\infty\ \mbox{and}\ \int_0^1{\frac{1}{W(x)+C+1}\mu(dx)}<+\infty.
 \end{equation*}
\end{hypothese}
\begin{hypothese}[H3]
\begin{equation*}
 \int_1^{+\infty}{e^{-2Q(x)}dx}<+\infty\ \mbox{and}\ \int_0^1{x\:e^{-Q(x)}dx}<+\infty.
\end{equation*}
\end{hypothese}

\begin{theoreme}
 \label{th1}
 Assume that hypotheses (H1) and (H2 or H3) are satisfied. Then
\begin{equation*}
 \nu_{\epsilon}\FlechBH{\longrightarrow}{\epsilon\rightarrow0}{}\nu\in{\cal M}_1(]0,+\infty[),
\end{equation*}
where $\nu$ is a QSD for $\mathbb{P}^0$, which is equal to the Yaglom limit $\lim_{t\rightarrow+\infty}{\mathbb{P}^0_x(X_t\in.|t<\tau_{\partial})}$, $\forall x\in]0,+\infty[$.
\end{theoreme}

\begin{remarque}
 \upshape
 The hypotheses (H1) and (H2 or H3) are the assumptions that are made in \cite{ca1} to prove the existence of the Yaglom limit.
\end{remarque}

\begin{remarque}
 \upshape
 If a process satisfies the hypotheses of Theorem \ref{th1}, then it is killed in finite time a.s or it is never killed a.s.
Indeed, assume that the process can be killed in finite time with a positive probability. Then $\int_1^{0}{e^{Q(x)}\left(\int_1^x{e^{-Q(y)}dy}\right)dx}<+\infty$ (see \cite{ik1}) and $\int_1^{+\infty}{e^{Q(x)}dx}=+\infty$ (as a consequence of (H1) and (H2 or H3)). But this two conditions are fulfilled if and only if the process is killed in finite time almost surely (see \cite[Theorem 3.2 p.450]{ik1}).
\end{remarque}

\begin{remarque}
\label{re1}
\upshape
The existence of a QSD for $\mathbb{P}^0$ can be seen as a consequence of Theorem \ref{th1}.
The existence of the Yaglom limit is proved in \cite[Theorem 5.2]{ca1}.
\end{remarque}

\begin{remarque}
 \upshape
In Part \ref{par3}, we give the counterpart of Theorem \ref{th1} for diffusions with values in a bounded interval.
\end{remarque}

The end of the section is devoted to the proof of Theorem \ref{th1}.

\subsection{Tightness of the family $(\nu_{\epsilon})_{0<\epsilon<1/2}$}

This part is devoted to the proof of the following result,
\begin{proposition}
\label{pr1}
Assume that the hypotheses (H1) and (H2 or H3) are satisfied. Then the family $(\nu_{\epsilon})_{0<\epsilon\leq1/2}$ is tight. Moreover, every limit point is absolutely continuous with respect to the Lebesgue measure.
\end{proposition}

We know that $L^{\epsilon}\eta_{\epsilon}=-\lambda_{\epsilon}\eta_{\epsilon}$, $\eta_{\epsilon}\in C^2([\epsilon,1/\epsilon])$, and $\eta_{\epsilon}$ satisfies the differential equation
\begin{equation*}
 1/2 \eta_{\epsilon}''(x)-q(x)\eta'_{\epsilon}(x)=-\lambda_{\epsilon} \eta_{\epsilon}(x)
\end{equation*}
with the boundary conditions
\begin{equation*}
 \eta_{\epsilon}(\epsilon)=\eta_{\epsilon}(1/\epsilon)=0.
\end{equation*}

Define $v_{\epsilon}=\eta_{\epsilon}e^{-Q}$. By \eqref{eq70}, we know that:
\begin{equation*}
 \int_{\epsilon}^{1/\epsilon}{v_{\epsilon}(x)^2 dx}=1.
\end{equation*}
We have
\begin{equation*}
 v_{\epsilon}(x)W(x)-v_{\epsilon}''(x)=2\lambda_{\epsilon} v_{\epsilon}(x),
\end{equation*}
with the boundary conditions
\begin{equation}
\label{eqc2}
 v_{\epsilon}(\epsilon)=v_{\epsilon}(1/\epsilon)=0.
\end{equation}

\begin{lemme}
\label{le5}
Assume that the  hypothesis (H1) is fulfilled.
Then $(v_{\epsilon})_{0<\epsilon<1/2}$ is uniformly bounded above and the family $(v_{\epsilon}^2(x) dx)_{0<\epsilon<1/2}$ is tight.
\end{lemme}
\textit{Proof of Lemma \ref{le5} : }
From the differential equation satisfied by $v_{\epsilon}$, we have
\begin{equation*}
 v_{\epsilon}(x)^2W(x)-v_{\epsilon}''(x)v_{\epsilon}(x)=2\lambda_{\epsilon} v_{\epsilon}(x)^2.
\end{equation*}
Integrating by parts and looking at the boundary conditions (\ref{eqc2}),
\begin{equation*}
 \int_{\epsilon}^{1/\epsilon}{v_{\epsilon}''(x)v_{\epsilon}(x)}dx=-\int_{\epsilon}^{1/\epsilon}{v_{\epsilon}'(x)^2}dx,
\end{equation*}
where $v_{\epsilon}$ is normalized in $L^2(dx)$. That implies
\begin{eqnarray*}
 \int_{\epsilon}^{1/\epsilon}{v_{\epsilon}'(x)^2}dx+\int_{\epsilon}^{1/\epsilon}{v_{\epsilon}(x)^2 W(x)}dx&=&2\lambda_{\epsilon}.
\end{eqnarray*}
The eigenvalue $\lambda_{\epsilon}$ of $-L^{\epsilon}$ is given by (see for instance \cite[chapter XI, part 8]{yo1})
\begin{eqnarray}
 \lambda_{\epsilon}&=&\inf_{\phi\in C^{\infty}_0\left(]\epsilon,1/\epsilon[\right)}{(L^{\epsilon}\phi,\phi)_{\mu}},\nonumber\\
                   &=&\inf_{\phi\in C^{\infty}_0\left(]\epsilon,1/\epsilon[\right)}{(L^0\phi,\phi)_{\mu}},\label{eq69}
\end{eqnarray}
where $C^{\infty}_0\left(]\epsilon,1/\epsilon[\right)$ is the vector space of infinitely differentiable functions with compact support in $]\epsilon,1/\epsilon[$ and  $(f,g)_{\mu}=\int_0^{+\infty}{f(u)g(u)d\mu(u)}$. We deduce from it that $\lambda_{\epsilon}$ increases with $\epsilon$ and is uniformly bounded above by $\lambda_{1/2}$.

We have then
\begin{equation}
\label{eqc3}
 0\leq \int_{\epsilon}^{1/\epsilon}{v_{\epsilon}'(x)^2}dx+\int_{\epsilon}^{1/\epsilon}{v_{\epsilon}(x)^2 (W(x)+C+1)}dx\leq 2\lambda_{1/2}+C+1.
\end{equation}
Looking at the boundary conditions (\ref{eqc2}), we obtain, for all $x\in]\epsilon,1/\epsilon[$,
\begin{eqnarray*}
 v_{\epsilon}^2(x)&=&-2\int_{x}^{1/\epsilon}{v_{\epsilon}'(y)v_{\epsilon}(y)}dy\\
       &\leq&-\frac{2}{\min_{[x,1/\epsilon[}{\sqrt{W+C+1}}}\int_x^{1/\epsilon}
             {v_{\epsilon}'(y)v_{\epsilon}(y)\sqrt{W(y)+C+1}}dy.
\end{eqnarray*}
Then, applying the Cauchy-Schwarz inequality to the right term above,
\begin{eqnarray*}
 v_{\epsilon}^2(x)&\leq&\frac{2}{\min_{[x,1/\epsilon[}{\sqrt{W+C+1}}}\sqrt{\int_{x}^{1/\epsilon}
             {v_{\epsilon}'(y)^2}dy}\sqrt{\int_{x}^{1/\epsilon}{v_{\epsilon}(y)^2(W(y)+C+1)}dy}.
\end{eqnarray*}
From (\ref{eqc3}), the integral product is bounded by $2\lambda_{1/2}+C+1$, thus $\exists A>0$, independent from $\epsilon$, such that
\begin{eqnarray}
\label{eqc4}
 v_{\epsilon}^2(x)&\leq&\frac{A}{\min_{[x,1/\epsilon[}{\sqrt{W+C+1}}}\nonumber\\
       &\leq&\frac{A}{\min_{[x,+\infty[}{\sqrt{W+C+1}}}\label{eqc7},
\end{eqnarray}
where $W(x)+C+1\geq 1$ for all $x\in]0,+\infty[$, thanks to Hypothesis (H1). That implies the first part of Lemma \ref{le5}.

Let us prove that the family $(v_{\epsilon}^2 dx)_{0<\epsilon<1/2}$ is tight. Fix $\delta>0$. We have to find  a compact subset $K_{\delta}$ in $]0,+\infty[$ such that 
\begin{equation}
 \label{eq67}
 \int_{]0,+\infty[\setminus K_{\delta}}{v_{\epsilon}^2(x)dx}\leq\delta,
\end{equation}
for all $\epsilon\in]0,1/2[$.
Thanks to \eqref{eqc7}, we have $v_{\epsilon}^2\leq A$, then
\begin{equation*}
 \int_0^{\delta/(2A)}{v_{\epsilon}^2}\leq \delta/2.
\end{equation*}
From the second part of Hypothesis (H1), $\exists M_{\delta}>0$ such that $W(x)+C+1>2(2\lambda+C+1)/\delta$ for all $x\geq M_{\delta}$. That implies
\begin{eqnarray*}
 \int_{M_{\delta}}^{+\infty}{v_{\epsilon}(x)^2}dx&\leq&\int_{\epsilon}^{1/\epsilon}{v_{\epsilon}(x)^2 \frac{\delta(W(x)+C+1)}{2(2\lambda+C+1)}}dx\\
                                                             &\leq&\delta/2,
\end{eqnarray*}
where the last inequality is due to \eqref{eqc3}.
Finally, the compact set $K_{\delta}=[\delta/(2A),M_{\delta}]$ satisfies \eqref{eq67}.
$\Box$

\begin{lemme}
\label{le6}
Assume that (H1) is satisfied.
Then $\int_{\epsilon}^{1/\epsilon}{\eta_{\epsilon}(y)d\mu(y)}$ is uniformly bounded below by a constant $B>0$.
\end{lemme}
\textit{Proof of Lemma \ref{le6} :} Assume that $\int_{\epsilon}^{1/\epsilon}{\eta_{\epsilon}(y)d\mu(y)}$ isn't uniformly bounded below : one can find a sub-sequence $\int_{\epsilon_k}^{1/\epsilon_k}{v_{\epsilon_k}(y)e^{-Q(y)}dy}=\int_{\epsilon_k}^{1/\epsilon_k}{\eta_{\epsilon_k}(y)d\mu(y)}$, where $\epsilon_k\rightarrow0$, which tends to $0$.
From Lemma \ref{le5}, $(v_{\epsilon})_{0<\epsilon\leq1/2}$ is uniformly bounded, so that
$\int_{\epsilon_k}^{1/\epsilon_k}{v_{\epsilon_k}(y)^2 e^{-Q(y)}dy}\rightarrow0$. The family $(v_{\epsilon}(x)^2dx)$ being tight, one can find (after extracting a sub-sequence) a positive map $m$ such that, for all continuous and bounded $\phi:\mathbb{R}_+\rightarrow\mathbb{R}$,
\begin{equation}
\label{eqc5}
 \int_{\epsilon_k}^{1/\epsilon_k}{v_{\epsilon_k}(y)^2 \phi(y)dy}\rightarrow \int_{0}^{+\infty}{m(y) \phi(y) dy}.
\end{equation}
Indeed, $(v_{\epsilon}^2)$ being uniformly bounded, all limit measure is absolutely continuous with respect to the Lebesgue measure. In particular,
\begin{equation*}
 \int_{\epsilon_k}^{1/\epsilon_k}{v_{\epsilon_k}(y)^2 \min{(e^{-Q(y)},1)}dy}\rightarrow \int_{0}^{+\infty}{m(y)\min{(e^{-Q(y)},1)}dy},
\end{equation*}
then
\begin{equation*}
 \int_{0}^{+\infty}{m(y)\min{(e^{-Q(y)},1)}dy}=0.
\end{equation*}
But $\min{(e^{-Q(.)},1)}$ is continuous and positive on $\mathbb{R}_+$, so that $m$ vanishes almost every where. Finally, by the convergence property (\ref{eqc5}) applied to $\phi$ equal to $1$ almost everywhere, we have
\begin{equation*}
 1=\int_{\epsilon_k}^{1/\epsilon_k}{v_{\epsilon_k}^2dx}\rightarrow 0,
\end{equation*}
what is absurd.
Thus, one can define $B=\inf_{\epsilon}{\int_{\epsilon}^{1/\epsilon}{\eta_{\epsilon}(y)d\mu(y)}}/A>0$. $\Box$

\begin{lemme}
 \label{le7}
 Assume that (H1) and (H2) are satisfied.
Then the family $(\eta_{\epsilon}(x)d\mu(x))_{0<\epsilon<1/2}$ is tight.
\end{lemme}
\textit{Proof of Lemma \ref{le7} :} By \eqref{eqc3}, we have
\begin{eqnarray*}
 \int_{\epsilon}^{1/\epsilon}{\eta_{\epsilon}^2(y)(W(y)+C+1)d\mu(y)}&=&\int_{\epsilon}^{1/\epsilon}{v_{\epsilon}^2(y)(W(y)+C+1)dy}\\
&\leq&\lambda_{1/2}+C+1.
\end{eqnarray*}
For all $\delta,M>0$, using Cauchy-Schwarz inequality, we get on one hand
\begin{eqnarray}
\label{eq71}
 \int_0^{\delta}{\eta_{\epsilon}(y)d\mu(y)}&\leq&\left(\int_0^{\delta}{\eta_{\epsilon}(y)^2(W(y)+C+1)d\mu(y)}\right)^{\frac{1}{2}}\left(\int_0^{\delta}{\frac{1}{W(y)+C+1}d\mu(y)}\right)^{\frac{1}{2}}\\
&\leq&\left(\lambda_{1/2}+C+1\right)^{\frac{1}{2}}\left(\int_0^{\delta}{\frac{1}{W(y)+C+1}d\mu(y)}\right)^{\frac{1}{2}}.
\end{eqnarray}
On the other hand,
\begin{eqnarray}
\label{eq72}
 \int_M^{+\infty}{\eta_{\epsilon}(y)d\mu(y)}&\leq&\left(\int_M^{+\infty}{\eta_{\epsilon}^2(y)d\mu(y)}\right)^{\frac{1}{2}}\left(\int_M^{+\infty}{d\mu(y)}\right)^{\frac{1}{2}}\\
&\leq&\left(\int_M^{+\infty}{d\mu(y)}\right)^{\frac{1}{2}}.
\end{eqnarray}
Thanks to (H2), both terms are going to $0$ uniformly in $\epsilon$, when $\delta$ and $M$ tend respectively to $0$ and $+\infty$. As a consequence, the family $(\eta_{\epsilon}(x)d\mu(x))_{0<\epsilon<1/2}$ is tight. $\Box$

\begin{lemme}
 \label{le8}
 Assume that (H1) and (H3) hold.
Then the family $(\eta_{\epsilon}(x)d\mu(x))_{0<\epsilon<1/2}$ is tight.
\end{lemme}
\textit{Proof of Lemma \ref{le8} :} The first part of the hypothesis (H3) is the same as (H2)'s one, then
\begin{eqnarray*}
 \int_M^{+\infty}{\eta_{\epsilon}(y)d\mu(y)}\rightarrow0
\end{eqnarray*}
when $M$ goes to infinity, uniformly in $\epsilon$.

Moreover, there exists a constant $K>0$ such that, for any $x\in]0,1]$ and any $\epsilon\in]0,1/2]$,
\begin{eqnarray*}
\eta_{\epsilon}(x)\leq K x e^{Q(x)}.
\end{eqnarray*}
This is a consequence of \cite[Proposition 4.3]{ca1} whose proof is still available under our settings.
This inequality allows us to conclude the proof of Lemma \ref{le8}. $\Box$

Thanks to equality (\ref{eqb1}) and Lemmas \ref{le6}, \ref{le7} and \ref{le8}, the first part of Proposition \ref{pr1} is proved. Moreover, $\nu_{\epsilon}$ has a density with respect to the Lebesgue measure which is bounded on every compact set, uniformly in $\epsilon>0$. Thus every limit point is absolutely continuous with respect to the Lebesgue measure.

\subsection{The limit points of the family $(\nu_{\epsilon})_{0<\epsilon<1/2}$}
\begin{proposition}
\label{pr11}
Assume that Hypotheses (H1) and (H2 or H3)  are fulfilled and let $\nu$ be a probability measure which is the limit of a sub-sequence $(\nu_{\epsilon_k})_{k\in\mathbb{N}}$, where $\epsilon_k\rightarrow0$ when $k\rightarrow\infty$. Then $\nu$ is a QSD with respect to $\mathbb{P}^0$.
\end{proposition}
\textit{Proof of Proposition \ref{pr11} : }From Proposition \ref{pr1}, the family $(\nu_{\epsilon})_{0<\epsilon<1/2}$ is tight. Let $\nu$ be a limit point of the family $(\nu_{\epsilon})_{0<\epsilon<1/2}$. There exists a sub-sequence $(\nu_{\epsilon_k})_k$ which converges to $\nu$, where $(\epsilon_k)_{k\in\mathbb{N}}$ is a decreasing sequence which tends to $0$.
We already know that $\nu$ is absolutely continuous with respect to the Lebesgue measure. That implies that, for all open intervals $D=]c,d[\subset\mathbb{R}_+$,
\begin{equation}
\label{eqc12} 
 \nu_{\epsilon_k}(D)\rightarrow\nu(D),
\end{equation}
and, for all bounded maps $\phi$ continuous on $\mathbb{R}_+$,
\begin{equation}
\label{eqc11}
 \int_{\mathbb{R}_+}{\phi(x)d\nu_{\epsilon_k}(x)}\rightarrow\int_{\mathbb{R}_+}{\phi(x)d\nu(x)}.
\end{equation}

Let $\nu_{t}$ (resp. $\nu_{\epsilon,t}$) be the distribution at time $t$ of a diffusion with law $\mathbb{P}^0_{\nu}$ (resp. $\mathbb{P}^{\epsilon}_{\nu_{\epsilon}}$), conditioned to be not killed until time $t$, that is 
\begin{equation*}
\nu_{t}(dx)=\mathbb{P}^0_{\nu}(\omega_t\in dx|\tau_{\partial}>t)
\end{equation*}
and
\begin{equation*}
\nu_{\epsilon,t}(dx)=\mathbb{P}^{\epsilon}_{\nu_{\epsilon}}(\omega_t\in dx|\tau_{\partial}>t)
\end{equation*}
The probability measure $\nu_{\epsilon}$ being a QSD for $\mathbb{P}^{\epsilon}$, we have $\nu_{\epsilon,t}=\nu_{\epsilon}$ for all $t>0$. We want to show that $\nu=\nu_{t}$ for all $t>0$, we have then to prove the following convergence result:
\begin{equation}
 \label{eq68}
 \forall t>0,\ \forall D=]c,d[\subset\mathbb{R}_+,\ \nu_{\epsilon_k,t}(D)\FlechBH{\longrightarrow}{k\rightarrow+\infty}{}\nu_{t}(D).
\end{equation}
Indeed, suppose that (\ref{eq68}) holds, then on the one hand, $\nu_{\epsilon_k}(D)=\nu_{\epsilon_k,t}(D)\rightarrow\nu_{t}(D)$. In the other hand $\nu_{\epsilon_k}(D)\rightarrow\nu(D)$. We have then $\nu_{t}(]c,d[)=\nu(]c,d[)$, $\forall ]c,d[\subset\mathbb{R}_+$ and this conclude the proof of Proposition \ref{pr11}.

Let us prove (\ref{eq68}). By definition,
\begin{equation*}
 \nu_{\epsilon_k,t}(D)=\frac{\int_{]\epsilon_k,1/\epsilon_k[}{\mathbb{P}^{\epsilon_k}_x(\omega_t\in D)d\nu_{\epsilon_k}(x)}}
    {\int_{]\epsilon_k,1/\epsilon_k[}{\mathbb{P}^{\epsilon_k}_x(\tau_{\partial}>t)d\nu_{\epsilon_k}(x)}}.
\end{equation*}
The numerator is equal to
\begin{multline*}
 \int_{]\epsilon_k,1/\epsilon_k[}{\mathbb{P}^{\epsilon_k}_x(\omega_t\in D)d\nu_{\epsilon_k}(x)}=\int_{]\epsilon_k,1/\epsilon_k[}{\mathbb{P}^0_x(\omega_t\in D)d\nu_{\epsilon_k}(x)}\\
  +\int_{]\epsilon_k,1/\epsilon_k[}{\left[\mathbb{P}^{\epsilon_k}_x(\omega_t\in D)-\mathbb{P}^0_x(\omega_t\in D)\right]d\nu_{\epsilon_k}(x)}
\end{multline*}
For all $t>0$, the map $x\mapsto \mathbb{P}^0_x(\omega_t\in D)$ is continuous and bounded, then, by the convergence property (\ref{eqc11}),
\begin{equation*}
 \int_{]\epsilon_k,1/\epsilon_k[}{\mathbb{P}^0_x(\omega_t\in D)d\nu_{\epsilon_k}(x)}\rightarrow\int_{]0,\infty[}{\mathbb{P}^0_x(\omega_t\in D)d\nu(x)},
\end{equation*}
Assume that (H2) is fulfilled, then, similarly to (\ref{eq71}) and (\ref{eq72}), we have, for all bounded continuous functions $f:]0,+\infty[\rightarrow\mathbb{R}$ and all $M>0$,
\begin{equation*}
 \int_M^{+\infty}{f(x)\eta_{\epsilon_k}(x)d\mu(x)}\leq\left(\int_M^{+\infty}{\left|f(x)\right|^2 d\mu(x)}\right)^{\frac{1}{2}},
\end{equation*}
and, for all $m>0$,
\begin{equation*}
 \int_0^m{f(x)\eta_{\epsilon_k}(x)d\mu(x)}\leq\left(\lambda_{1/2}+C+1\right)^{\frac{1}{2}}\left(\int_0^m{\frac{\left|f(x)\right|^2}{W(x)+C+1}d\mu(x)}\right)^\frac{1}{2}.
\end{equation*}
Replacing $f(x)$ by $\mathbb{P}^0_x(\omega_t\in D)-\mathbb{P}^{\epsilon_k}_x(\omega_t\in D)$, which is decreasing to $0$ when $k\rightarrow\infty$, and by monotone convergence theorem, we have
\begin{equation*}
 \int_{]0,m[}{\left[\mathbb{P}^{\epsilon_k}_x(\omega_t\in D)-\mathbb{P}^0_x(\omega_t\in D)\right]d\nu_{\epsilon_k}(x)}\FlechBH{\longrightarrow}{k\rightarrow+\infty}{}0
\end{equation*}
and
\begin{equation*}
 \int_{]M,+\infty[}{\left[\mathbb{P}^{\epsilon_k}_x(\omega_t\in D)-\mathbb{P}^0_x(\omega_t\in D)\right]d\nu_{\epsilon_k}(x)}\FlechBH{\longrightarrow}{k\rightarrow+\infty}{} 0.
\end{equation*}
Finally, the density of $\nu_{\epsilon_k}$ being bounded above in every compact set $[m,M]$, uniformly in $\epsilon_k$, the same argument of monotone convergence gives us
\begin{equation*}
 \int_{]\epsilon_k,1/\epsilon_k[}{\left[\mathbb{P}^{\epsilon_k}_x(\omega_t\in D)-\mathbb{P}^0_x(\omega_t\in D)\right]d\nu_{\epsilon_k}(x)}\FlechBH{\longrightarrow}{k\rightarrow+\infty}{} 0.
\end{equation*}
With similar arguments, the same holds under (H3). Finally, we obtain
\begin{equation*}
 \int_{]\epsilon_k,1/\epsilon_k[}{\mathbb{P}^{\epsilon_k}_x(\omega_t\in D)d\nu_{\epsilon_k}(x)}\rightarrow\int_{]0,+\infty[}{\mathbb{P}^0_x(\omega_t\in D)d\nu(x)}.
\end{equation*}

Thanks to \cite[Lemma 5.3 and Theorem 2.3]{ca1}, the map $x\mapsto \mathbb{P}^0_x(\tau_{\partial}>t)=\mathbb{P}^0_x(\omega_t\in]0,+\infty[)$ is continuous, and $\mathbb{P}^0_x(\tau_{\partial}>t)-\mathbb{P}_x^{\epsilon_k}(\tau_{\partial}>t)$ is increasing to $0$ when $k\rightarrow\infty$. Thus the denominator can be treated in the same way. $\Box$

We can now conclude the proof of Theorem \ref{th1}:
\begin{proposition}
 \label{pr12}
 Assume that (H1) and (H2 or H3) hold.
 The limit measure $\nu$ in the statement of Proposition \ref{pr11} is unique. Moreover $\nu$ is the Yaglom limit associated with $\mathbb{P}^0_x$, $\forall x\in]0,+\infty[$.
\end{proposition}
\textit{Proof of Proposition \ref{pr12} : }The proof of Proposition \ref{pr11} implies that
\begin{equation*}
 \mathbb{P}_{\nu_{\epsilon}}^{\epsilon}(\tau_{\partial}>t)\FlechBH{\longrightarrow}{\epsilon\rightarrow0}{}\mathbb{P}_{\nu}^{0}(\tau_{\partial}>t),\ \forall t>0.
\end{equation*}
The probability measure $\nu_{\epsilon}$ being a QSD for $\mathbb{P}^{\epsilon}$, we have
\begin{equation*}
 \mathbb{P}_{\nu_{\epsilon}}^{\epsilon}(\tau_{\partial}>t)=e^{-\lambda_{\epsilon}t},\ \forall t>0.
\end{equation*}
Thanks to (\ref{eq69}), $\lambda_{\epsilon}$ is decreasing to $\lambda_0=\inf_{\phi\in C^{\infty}_0\left(]0,+\infty[\right)}{(L^0\phi,\phi)_{\mu}}$ when $\epsilon$ goes to $0$. As a consequence,
\begin{equation*}
 \mathbb{P}_{\nu}^{0}(\tau_{\partial}>t)=e^{-\lambda_0 t},\ \forall t>0.
\end{equation*}
In this case, the density of $\nu$ with respect to $d\mu$ is an eigenfunction of $L^0$ with eigenvalue $-\lambda_0<0$, where $(L^0)^*$ is the adjoint operator of $L^0$ (this is a consequence of the spectral decomposition proved in \cite[Theorem 3.2]{ca1}). As defined, $-\lambda_0$ is at the bottom of the spectrum of $(L^0)^*$. Thanks to \cite[Theorem 3.2]{ca1}, this eigenvalue is simple. Moreover, \cite[Theorem 5.2]{ca1} states that this QSD is equal to $\lim_{t\rightarrow+\infty}{\mathbb{P}^0_x(X_t\in.|\tau_{\partial}>t)}$, what concludes the proof. $\Box$

\subsection{Diffusions with values in a bounded interval}
\label{par3}
Theorem \ref{th1} is stated for $1$-dimensional diffusions with values in $]0,+\infty[$. However, most of the proofs can be easily adapted to diffusions with values in a bounded interval $]a,b[$, where $-\infty<a<b<+\infty$, defined by the SDE
\begin{equation*}
\label{eq76}
 dX_t=dB_t-q(X_t)dt,\ X_0=x\in]a,b[,
\end{equation*}
and killed when it hits $a$ or $b$.
Here $B$ is a standard Brownian motion and $q\in C^1(]a,b[)$.

More precisely, let us denote by $\mathbb{P}^0$ the law of such a diffusion. For each $\epsilon>0$, define $\mathbb{P}^{\epsilon}$ as the law of a diffusion with values in $]a+\epsilon,b-\epsilon[$, driven by the SDE
\begin{equation*}
 dX_t=dB_t-q(X_t)dt,\ X_0=x\in]a+\epsilon,b-\epsilon[,
\end{equation*}
and killed when it hits $a+\epsilon$ or $b-\epsilon$. As proved in \cite{pi1}, there exists a unique QSD $\nu_{\epsilon}$ associated with $\mathbb{P}^{\epsilon}$.

We define $Q(x)=\int_{(a+b)/2}^{x}{q(y)dy}$ and $W(x)=q(x)^2-q'(x)$. The counterpart of Theorem \ref{th1} under these settings is
\begin{theoreme}
\label{th10}
 Assume that the following hypotheses are fulfilled:
\begin{hypothese}[HH1]
 $W$ is uniformly bounded below by $-C$, where $C$ is a positive constant.
\end{hypothese}
\begin{hypothese}[HH2]
 $x\mapsto \frac{1}{W(x)+C+1}e^{-2Q(x)}$ or $x\mapsto (x-a)e^{-Q(x)}$ is integrable on a neighbourhood of $a$.
\end{hypothese}
\begin{hypothese}[HH3]
 $x\mapsto \frac{1}{W(x)+C+1}e^{-2Q(x)}$ or $x\mapsto (b-x)e^{-Q(x)}$ is integrable on a neighbourhood of $b$.
\end{hypothese}
Then the family of QSD $(\nu_{\epsilon})$ is tight as family of measures on $]a,b[$.
Moreover, every limit point of the family $(\nu_{\epsilon})_{0<\epsilon<1/2}$ is a QSD for $\mathbb{P}^0$.
\end{theoreme}

\begin{remarque}
 \upshape
Our aim isn't to develop this part, but we point out that to show that Proposition \ref{pr12} remains valid, most of the arguments used to prove the key results \cite[Theorem 3.2]{ca1} and \cite[Theorem 5.2]{ca1} can be adapted to these settings.
\end{remarque}

\section{Approximation of $\nu_{\epsilon}$, QSD for $\mathbb{P}^{\epsilon}$}
\label{par2}
We are interested in proving an approximation method for the QSD associated with $\mathbb{P}^{\epsilon}$. It will be sufficient to prove it for
any diffusion $(X_t)$ taking its values in $]0,1[$ in place of $]\epsilon,1/\epsilon[$, defined by the stochastic differential equation (SDE)
\begin{eqnarray}
 \label{eq1}
 dX_t=dB_t-q(X_t)dt,\ X_0=x>0,
\end{eqnarray}
and killed when $X_t$ hits the boundary $\{0,1\}$. Here $B$ is a real Brownian motion and $q\in C^1([0,1])$. The law of $X$ will be denoted by $\mathbb{P}$. 

From \cite{pi1}, the QSD of $X$ is unique and equals the Yaglom limit. It will be denoted by $\nu$. For notational convenience, new notations have been defined for this section, which is totally independent of the previous one.

Fix $N\geq2$ and let us define formally the interacting particle process with $N$ particles described in the introduction. Let $B^1,...,B^N$ be  $N$ independent Brownian motions and $(X^1_0,...,X^N_0)\in]0,1[^N$ be the starting point of the process.
\begin{itemize}
 \item For each $i\in\{1,...,N\}$, the particle $X^i$ evolves in $]0,1[$ and satisfies the SDE $dX^i_t=dB^i_t-q(X^i_t)dt$ (and then it is independent of the others) until $\tau^i_1=\inf\{t\geq0,\ X^i_{t-}=0\ \mbox{or}\ 1\}$.
 \item At time $\tau_1=\min\{\tau^1_1,...,\tau^1_N\}$, the path of a particle, denoted by $i_1$ (it is unique), has a left limit equal to $0$ or $1$.
 \item A particle $j_1$ is chosen in $\{1,...,N\}\setminus\{i_1\}$. The particle $i_1$ jumps on the position of the particle $j_1$: we set $X^{i_1}_{\tau_1}:=X^{j_1}_{\tau_1}$.
 \item After time $\tau_1$, each particle $X^i$ evolves in $]0,1[$ with respect to the SDE $dX^i_t=dB^i_t-q(X^i_t)dt$ until $\tau^i_2=\inf\{t>\tau_1,\ X^i_{t-}=0\ \mbox{or}\ 1\}$. At time $\tau_1$, all the particles are in $]0,1[$, so that we have $\tau^i_2>\tau_1$ for all $i\in\{1,...,N\}$ almost surely.
 \item  At time $\tau_2=\min\{\tau^1_2,...,\tau^N_2\}$ (which is then strictly bigger than $\tau_1$), a unique particle $i_2$ has a path whose left limit is equal to $0$ or $1$.
 \item A particle $j_2$ is chosen in $\{1,...,N\}\setminus\{i_2\}$. The particle $i_2$ jumps on the position of the particle $j_2$: we set $X^{i_2}_{\tau_2}:=X^{j_2}_{\tau_2}$.
 \item After time $\tau_2$, the particles evolve independently from each other and so on.
\end{itemize}
Following this way, we define the strictly increasing sequence of stopping times $0<\tau_1<\tau_2<\tau_3<...$, the time $\tau_{\infty}=\lim_{n\rightarrow\infty}{\tau_n}$ and the interacting particle system $(X^1_t,...,X^N_t)$ for all $t\in[0,\tau_{\infty}[$. The law of $(X^1,...,X^N)$ will be denoted by $\mathbb{P}^{ipp}$.

\paragraph{}
We can now state the main result of this section:

\begin{theoreme}
\label{th2}
$(X^1,...,X^N)$ is well defined, that means $\tau_{\infty}=+\infty$ almost surely. It is geometrically ergodic, with unique stationary distribution $M^N$.

 Let ${\cal X}^N$ be the empirical stationary measure of the interacting particle process with $N$ particles, that is the empirical measure of a random vector $(x^1,...,x^N)\in]0,1[^N$ distributed with respect to the stationary measure $M^N$ of the process $(X^1,...,X^N)$. The sequence of random measures $({\cal X}^N)_{N\geq2}$ converges in law to the deterministic measure $\nu$, QSD of the process $X$.
\end{theoreme}

Subsection \ref{parts21} is devoted to prove a coupling which ensures the non-degeneracy of the particles at the boundary. A consequence will be that $\tau_{\infty}=+\infty$ almost surely. In Subsection \ref{parts22}, the process is studied in finite time. We prove that the empirical measure of the process $(X^1,...,X^N)$ at time $t$ converges, when $N$ goes to infinity, to the distribution of $X_t$ conditioned to not be killed until time $t$, which is $\mathbb{P}_{\mu_0}(X_t\in.|X_t\neq\partial)$ ($\mu_0$ denotes the limit of the empirical measure of the process at time $0$). We prove in Subsection \ref{parts23} that the interacting particle system $(X^1,...,X^N)$ is geometrically ergodic. We conclude by showing the convergence of the empirical stationary measure to the QSD $\nu$. 

\subsection{Existence and non-degeneracy at the boundary}
\label{parts21}
One of the most important fact when studying the interacting particle system is that the particles don't degenerate at the boundary. This is an evidence if the particles are independent. In our case, we will prove a coupling between $(X^1,...,X^N)$ and an other process $(Y^1,...,Y^N)$ whose components are independent identically distributed and don't degenerate at the boundary and such that, for all $i\in\{1,...,N\}$,
\begin{eqnarray}
\label{eqco1}
0\leq Y^i\leq X^i\leq 1-Y^i\leq 1\ \ \mbox{a.s.}
\end{eqnarray}
With this construction, the process $(X^1,...,X^N)$ doesn't degenerate at the boundary, because each of its particles $X^i$ is contained in $[Y^i,1-Y^i]$. This coupling will be useful in each step of the proof.

\subsubsection{Coupling's construction}
Define $Q=\sup_{x\in]0,1[}{|q(x)|}$ (we have $Q<+\infty$ by hypothesis) and fix $i\in\{1,...,N\}$. The process $Y^i$ is defined with values in $[0,1/3]$ by the SDE 
\begin{eqnarray*}
dY^i_t=dW^i_t-Qdt,\ Y^i(0)=\min{\left\{X^i(0),1-X^i(0),1/3\right\}},
\end{eqnarray*}
 with $0$ and $1/3$ as reflecting boundary (see \cite{ka1} for the definition of a reflected diffusion). The coupling inequality \eqref{eqco1} is fulfilled at time $t=0$. The Brownian motion $W^i$ will depend on $B^i$ and on the position of $X^i$.

If $X^i_t$ belongs to $[0,1/3]$, then $X^i_t\leq1/3$ and the second part of the coupling inequality is satisfied, independently of the choice of $W^i$. We only need to ensure that $X^i_t$ stays bigger than $Y^i_t$.  If $X^i_t$ belongs to $[2/3,1[$, we only need to ensure that it is smaller than $1-Y^i_t$. If it belongs to $]1/3,2/3[$, the coupling inequality is obviously fulfilled, thanks to the reflection of $Y^i$ on $1/3$.

\begin{figure}[htbp]
\begin{center}
\input{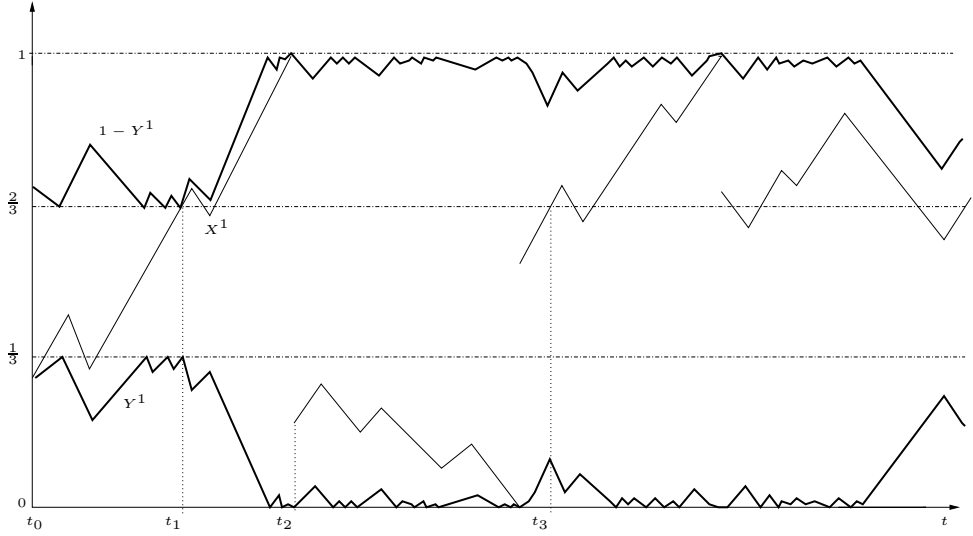}
\caption{The process $Y^1$}
\label{fig2}
\end{center}
\end{figure}

Assume that $X^i_t$ is in $]0,1/3]$ and that the coupling inequality is fulfilled at time $t$. We have $d(X^i_t-Y^i_t)=dB^i_t-dW^i_t+(Q-q(X_t))dt$, with $Q-q(X_t)\geq 0$. If we choose $W^i$ so that $dB^i_t-dW^i_t=0$, then $d(X^i_t-Y^i_t)=(Q-q(X_t))dt$ is increasing with time almost surely and the coupling inequality remains fulfilled. In a similar way, if $X^i_t$ is in $[2/3,1[$, we have to choose $W^i$ so that $dB^i_t+dW^i_t=0$. We will see in the proof of the coupling inequality that the jumps of $X^i$ and the reflexion of $Y^i$ on $0$ do not play any role in the coupling inequality (see Figure \ref{fig2}).

Let us define the Brownian motion $W^i$ in an explicit form. First, we build a sequence of strictly increasing times $(t^i_n)_{n\geq0}$ such that, for each $n\geq0$, $X^i_t\in]0,2/3[$ for all $t\in[t^i_{2n},t^i_{2n+1}[$ and $X^i_t\in]1/3,1[$ for all $t\in[t^i_{2n+1},t^i_{2n+2}[$. Define
\begin{align*}
	t^i_0=\inf{\{t\in[0,+\infty],\ X^i_t\in]0,1/3]\}},\\
	t^i_1=\inf{\{t\in[t_0,+\infty],\ X^i_t\in[2/3,1[\}},
\end{align*}
and, for $n\geq 1$,
\begin{align*}
	t^i_{2n}=\inf{\{t\in[t^i_{2n-1},+\infty],\ X^i_t\in]0,1/3]\}},\\
	t^i_{2n+1}=\inf{\{t\in[t^i_{2n},+\infty],\ X^i_t\in]2/3,1]\}}.
\end{align*}
Depending on the position of $X^i$, which is stated by the sequence $(t_n)_{n\geq2}$, we define $W^i$ by
\begin{eqnarray*}
	\label{eq8}
	W^i_t=-B^i_t\ \mbox{for}\ t\in[0,t_0],
\end{eqnarray*}
and, for all $n\geq0$,
\begin{eqnarray*}
	\label{eq9}
	W^i_t=W^i_{t_{2n}}+(B^i_t-B^i_{t_{2n}})\ \mbox{for}\ t\in[t_{2n},t_{2n+1}]\\
	\label{eq10}
	W^i_t=W^i_{t_{2n+1}}-(B^i_t-B^i_{t_{2n+1}})\ \mbox{for}\ t\in[t_{2n+1},t_{2n+2}].
\end{eqnarray*}
If $\lim_{n\rightarrow\infty}{t_n}<+\infty$, then $W^i$ can be extended by continuity on $[0,+\infty[$, so that 
$Y^i_t$ is well defined for all $t\in[0,+\infty[$.

The sequence $(t^i_n)_{n\in\mathbb{N}}$ is a sequence of stopping times for the natural filtration of the process $(X^1,...,X^N)$. Conditionally to the sequence $(j_n)_{n\geq1}$ (the particles which are chosen at each successive jump), $(X^1,...,X^N)$ only depends on $(B^1,...,B^N)$. One can then apply the strong Markov property to $(B^1,...,B^N)$ at time $t^i_n$, $\forall i\in\{1,...,N\}$ and $\forall n\in\mathbb{N}$.

As a direct consequence of the symmetry of a Brownian motion's law and of the strong Markov property applied to the $N$-dimensional Brownian motion $(B^1,...,B^i,...,B^N)$, $W^i$ is a Brownian motion. Note that the sequence of stopping times depends clearly on the position of other particles and then on the $B^j$'s, for $j\in \{1,...,N\}$. However, the Brownian motions $\{W^j\}_{j=1,...,N}$ are independent processes. Indeed, the sequence $(t_n)_{n\geq0}$ doesn't play any role in the law of $W^i$ and, to be convinced of that, one can compute the covariance matrix of the Gaussian variables $W^j_t$, for any $t\in[0,+\infty[$, which is clearly a diagonal one.

\begin{proposition}
\label{pr2}
The coupling inequality \eqref{eqco1} is fulfilled for all $t\in[0,\tau_{\infty}[$.
\end{proposition}
\textit{Proof of Proposition \ref{pr2} :} Define the time $\zeta=\inf\left\{0\leq t\leq\tau_{\infty}, Y^i_t > X^i_t\right\}$ and let us work conditionally to $\zeta<\tau_{\infty}$, then we have, by right continuity of the two processes, $Y^i_{\zeta}\geq X^i_{\zeta}$ a.s.

We first show that $\zeta$ is a jump time for the particle $i$. Assume the converse. If $\zeta=0$, then $X^i_{\zeta}=Y^i_{\zeta}$ and, if $\zeta>0$, then $X^i$ and $Y^i$ are continuous in a neighbourhood $[\zeta-h,\zeta+h]$ of $\zeta$, where $h>0$, and we have $Y^i_{t}\leq X^i_{t}$, for all $t\in[\zeta-h,\zeta[$, so that, by left continuity $Y^i_{\zeta}\leq X^i_{\zeta}$, and then $Y^i_{\zeta}=X^i_{\zeta}$. Therefore, we have $X^i_{\zeta}\in]0,1/3]$ and $\exists n\geq0$ such that $\zeta\in[t_{2n},t_{2n+1}[$, and, for $h$ small enough, one can assume that $\zeta+h<t_{2n+1}$. Then, for all $t\in[\zeta,\zeta+h]$, $d(X^i_t-Y^i_t)=(Q-q(X^i_t))dt+dL^{1/3}_t$, where $Q-q(X^i_t)\geq0$ and $L^{1/3}_t$ is an increasing process due to the reflecting property of the boundary $1/3$ for $Y^i$ (note that $Y^i_{\zeta}>0$, so that, for $h>0$ small enough, $Y^i_t>0$ for all $t\in[\zeta,\zeta+h[$). Then $X^i-Y^i$ stays non-negative between times $\zeta$ and $\zeta+h$, what contradicts the definition $\zeta$.

 The time $\zeta$ is then a jump time for the particle $i$. If $X^i_{\zeta-}=0$, then, by definition of $\zeta$, $Y^i_{\zeta-}=0$ and, by left continuity of the process $Y^i$, $Y^i_{\zeta}=0$, so that $X^i_{\zeta}=0$, what is impossible. Therefore $X^i_{\zeta-}=1>1-Y^i_{\zeta-}$, and, by existence of left limits for the two processes, $\exists t<\zeta$ such that $X^i_t>1-Y^i_t$. Define $\zeta'=\inf\left\{t\geq0, 1-Y^i_t < X^i_t\right\}$. We have then, conditionally to $\zeta<\tau_{\infty}$, $\zeta'<\zeta$.

By symmetry, conditionally to the event $\zeta'<\tau_{\infty}$, we have $\zeta<\zeta'$, then $\zeta<\tau_{\infty}$ and therefore $\zeta<\zeta'$. Finally, $\zeta=\zeta'=\tau_{\infty}$ almost surely.
 $\Box$

\subsubsection{Existence of the interacting particle process}

\begin{proposition}
\label{pr3}
 For all $N\geq2$, the interacting particle system $(X^1,...,X^N)$ is well defined, that is $\tau_{\infty}=+\infty$ almost surely.
\end{proposition}
\textit{Proof of Proposition \ref{pr3} :}
Let $N\geq2$ be the size of the interacting particle system and fix arbitrarily its starting point $x\in]0,1[^N$. We define the event $C_x=\{\tau_{\infty}<+\infty\}$.

Conditionally to $C_x$, the total number of jumps is equal to $+\infty$. There is a finite number of particles, then at least one particle, say $i_0$, makes an infinite number of jump before $\tau_{\infty}$. At each jump of $i_0$, a particle is uniformly chosen in $\{1,...,N\}$. By the law of large numbers, each particle is chosen infinitely often before $\tau_{\infty}$.
Assume that a particle, say $j_0$, remains  all the time in $]\epsilon,1-\epsilon[$, with $\epsilon>0$. $i_0$ will jump on the position of $j_0$ infinitely often. Then it will come back from $]\epsilon,1-\epsilon[$ to the boundary infinitely often in finite time, what is impossible. We deduce that, conditionally to $C_x$, all particles of the interacting particle system are going near to the boundary, that is
\begin{equation}
 \label{eq59}
 C_x\subset\left\{\lim_{t\rightarrow\tau_{\infty}}{\min{(X^i_t,1-X^i_t)}}=0\right\},
\end{equation}
for each particle $i\in\{1,2,...,N\}$.

Using the coupling inequality of Proposition \ref{pr2}, we deduce from (\ref{eq59}) that
\begin{equation*}
 \label{eq60}
 C_x\subset\left\{\lim_{t\rightarrow\tau_{\infty}}{(Y^1_t,...,Y^N_t)}=0\right\}.
\end{equation*}
Then, conditionally to $C_x$, $Y^1$ and $Y^2$ are independent reflected diffusions with bounded drifts, which hit $0$ at the same time. This occurs for two independent reflected Brownian motions with probability $0$, and then for $Y^1$ and $Y^2$ too, by the Girsanov's Theorem. That implies $P(C_x)=0$.

Finally, we have $\tau_{\infty}=+\infty$ almost surely.

\begin{remarque}\upshape
 One could hope to apply this method directly to the process with law $\mathbb{P}^0$ studied in the first part of this paper. Unfortunately, it can be very difficult to show the existence of the process or the non-degeneracy at the boundary: the drift being unbounded, the law of the reflected diffusion used in this proof isn't absolutely continuous with respect to the law of the reflected Brownian motion and stay at $0$ all the time after hitting it.
\end{remarque}

\subsubsection{Non-degeneracy at the boundary}
For all $r>0$, we define the open set $D_r=]r,1-r[$. Let $\mu^N(t,dx)$ (resp. $\mu'^N(t,dx)$) be the empirical measure of the system of particles $(X^i_t)_{i=1,...,N}$ (resp. $(Y_t^i)_{i=1,...,N}$), that is
\begin{eqnarray*}
 \mu^N(t,dx)=\frac{1}{N}\sum_{i=1}^N{\delta_{X^i_t}(dx)}\ \mbox{and}\ \mu'^N(t,dx)=\frac{1}{N}\sum_{i=1}^N{\delta_{Y^i_t}(dx)}.
\end{eqnarray*}
We will suppose that, at time $0$, the sequence of empirical measures $(\mu^N(0,dx))_{N\geq2}$ satisfies the following non-degeneracy property, which ensures that the mass of $\mu^N(0,dx)$ doesn't degenerate at the boundary, uniformly in $N$:
\begin{definition}
 The family of random probabilities $\{\mu^N(dx)\}$ is said to verify the non-degeneracy property if, for any $\epsilon>0$,
\begin{equation}
 \label{eq62}
 \lim_{r\rightarrow0}{\limsup_{N\rightarrow\infty}{P\left(\mu^N(D_r^c)>\epsilon\right)}}=0,
\end{equation}
where $D_r^c=]0,r]\cup[1-r,1[$.
\end{definition}
From definition of $(Y^i_0)_{i\in\{1,...,N\}}$, the non-degeneracy of $\mu'^N(0,dx)$ is the consequence of the non-degeneracy of $\mu^N(0,dx)$.
The end of the section is devoted to prove the following Proposition, which states that the non-degeneracy property is maintained over time:
\begin{proposition}
\label{pr4}
Assume that $(\mu^N(0,dx))_{N\geq 2}$ satisfies the non-degeneracy property, then, for all $T>0$ and all $\epsilon>0$,
\begin{eqnarray*}
	\label{eq11}
	\lim_{r\rightarrow 0}{\limsup_{N\rightarrow\infty}{P(\sup_{t\in[0,T]}{\mu^N(t,D_r^c)}>\epsilon)}}=0.
\end{eqnarray*} 
\end{proposition}
\textit{Proof of Proposition \ref{pr4} :} Fix $T>0$ and $\epsilon>0$.
Because of the non-degeneracy of $\mu'^N(0,dx)$, one can find $a>0$ such that
\begin{equation}
 \label{eq65}
 P\left(\frac{1}{N}\sum_{i=1}^N{\mathbf{1}_{Y_0^{i}\leq a}}\geq \epsilon/2\right)\rightarrow0
\end{equation}
when $N$ goes to $\infty$. 

We want to apply a law of large numbers, but the $Y_0^i$ aren't independent and depend on the number of particles $N$.
Let us define the diffusion $Z^i$ with values in $[0,1/3]$, defined by the SDE
\begin{eqnarray*}
dZ^i_t=dW'^i_t-Qdt,\ Z^i(0)=a,
\end{eqnarray*}
with $0$ and $1/3$ as reflecting boundary. Here the $W'^i$ are independent Brownian motions. The random processes $\mathbf{1}_{Z^i_.\in D_r^c}$ are independent, identically distributed with values in $D([0,T],\mathbb{R})$. One can apply to them the following law of large numbers, proved in \cite{ra1}:
\begin{equation*}
 \sup_{t\in[0,T]}\left|\frac{1}{N}\sum_{i=1}^N{\mathbf{1}_{Z^i_t\in D_r^c}}-E\left(\mathbf{1}_{Z^i_t\in D_r^c}\right)\right| \FlechBH{\longrightarrow}{N\rightarrow\infty}{Prob} 0.
\end{equation*}
We have $E\left(\mathbf{1}_{Z^i_.\in D_r^c}\right)=P\left(Z^i_.\in D_r^c\right)$, which tends uniformly to $0$ when $r\rightarrow 0$. We deduce from it that
\begin{equation}
 \label{eq66}
 \lim_{r\rightarrow0}\lim_{N\rightarrow\infty}P\left(\sup_{t\in[0,T]}{\frac{1}{N}\sum_{i=1}^N{\mathbf{1}_{Z^i_t\in D_r^c}}}\geq \epsilon/2\right)=0.
\end{equation}

For each number of particles $N\geq2$, one can easily find a coupling between $Z^i$ and $Y^i$, where $Z^i_0\leq Y^i_0$ implies $Z^i_t\leq Y^i_t$ for all $t\in[0,T]$. With such a coupling, we have
\begin{equation*}
\frac{1}{N}\sum_{i=1}^N{\mathbf{1}_{Y^i_t\in D_r^c}\mathbf{1}_{Y^i_0\geq a}} \leq \frac{1}{N}\sum_{i=1}^N{\mathbf{1}_{Z^i_t\in D_r^c}}.
\end{equation*}
By adding the contribution of the $Y^i$ which starts in $]0,a[$, we get
\begin{equation*}
\frac{1}{N}\sum_{i=1}^N{\mathbf{1}_{Y^i_t\in D_r^c}} \leq \frac{1}{N}\sum_{i=1}^N{\mathbf{1}_{Z^i_t\in D_r^c}}+\frac{1}{N}\sum_{i=1}^N{\mathbf{1}_{Y^i_0<a}}.
\end{equation*}
The limits \eqref{eq65} and \eqref{eq66} allow us to conclude the proof.

\subsection{Convergence in finite time}
\label{parts22}
Fix $T>0$. This section is devoted to the proof of the following proposition, which states that the empirical measure process converges to the distribution of the process $X$ conditioned to not be killed.
\begin{proposition}
\label{pr5}
 Assume that $(\mu^N(0,dx))_{N\in\mathbb{N}}$ converges in law to the random probability measure $\mu(0,dx)$ with respect to the weak topology and satisfies the non-degeneracy property.

Then, $\forall T>0$, the measure processes $(\mu^N(t,dx))_{t\in[0,T]}$ converge in law to $(\mathbb{P}(X_t\in dx|X_t\neq\partial))_{t\in[0,T]}$ in the Skorokhod space $D([0,T],{\cal M}_1(]0,1[))$ when $N\rightarrow\infty$. Here ${\cal M}_1(]0,1[)$ denotes the space of probability measures on $]0,1[$ equipped with the weak topology.
\end{proposition}
\textit{Proof of Proposition \ref{pr5} :}
For all maps $\psi\in C^2_0([0,1])$ vanishing on $\{0,1\}$, one can apply the It\^o's formula to the semimartingale 
$\psi(X^i_t)$ (see \cite[Theorem 27.1]{me2}), whose number of jumps in $[0,T]$ is finite almost surely :
\begin{multline}
 \psi(X^i_t)=\psi(X^i_0)+\int_0^t{\psi'(X^i_{s\minus}) dB^i_s}+\int_0^t{\left(\psi'(X^i_{s\minus})q(X^i_{s\minus})+\frac{1}{2}\psi''(X^i_{s\minus})\right) ds}\\
+\sum_{0\leq s\leq t}{\psi(X^i_s)-\psi(X^i_{s\minus})}
\end{multline}
Let us denote by $(\tau^i_n)_{1\leq n}$ the increasing sequence of jump times of the particle $i$. We have
\begin{align}
 \sum_{0\leq s\leq t}{\psi(X^i_s)-\psi(X^i_{s\minus})}=&\sum_{0\leq \tau^i_n\leq t}{\psi(X^i_{\tau_n})-\psi(X^i_{\tau_n\minus})}\\
                                                 =&\sum_{0\leq \tau^i_n\leq t}{\psi(X^i_{\tau_n})},
\end{align}
because $X^i_{\tau^i_n\minus}\in\{0,1\}$ and $\psi(0)=\psi(1)=0$. That implies
\begin{align}
\sum_{0\leq s\leq t}{\psi(X^i_s)-\psi(X^i_{s\minus})}=&\sum_{0\leq \tau^i_n\leq t}{\left(\psi(X^i_{\tau^i_n})-\frac{1}{N-1}\sum_{j=1}^N{\psi(X^j_{\tau^i_n\minus})}\right)}\\
 &+\frac{1}{N-1}\sum_{j=1}^N{\sum_{0\leq \tau^i_n\leq t}{\psi(X^j_{\tau^i_n\minus})}}.
\end{align}
By summing over $i\in\{1,...,N\}$, we obtain
\begin{multline}
 \label{eqd5}
 \mu^N(t,\psi)=\mu^N(0,\psi)+\int_0^t{\mu^N(s\minus,\psi'q+\frac{1}{2}\psi'')ds}+{\cal M}^c(\psi,t)+{\cal M}^j(\psi,t)\\+\frac{1}{N-1}\sum_{0\leq \tau_n\leq t}{\mu^N(\tau_n\minus,\psi)},
\end{multline}
where ${\cal M}^c(\psi,t)$ is the continuous martingale $\frac{1}{N}\sum_{i=1}^N{\int_0^t{\psi'(X^i_{s\minus})}dB^i_s}$ and ${\cal M}^j(\psi,t)$ is the pure jump martingale
\begin{equation}
 {\cal M}^j(\psi,t)=\frac{1}{N}\sum_{i=1}^N{\sum_{0\leq \tau^i_n\leq t}{\left(\psi(X^i_{\tau^i_n})-\frac{1}{N-1}\sum_{j=1}^N{\psi(X^j_{\tau^i_n\minus})}\right)}}.
\end{equation}

Now, we interpret each jump as a killing. Then we introduce a loss of $1/N$ of the total mass at each jump: we look at 
the measure process $\mu^N$ decreased by a factor $\frac{N-1}{N}$ at each jump.
More precisely, we set 
\begin{eqnarray*}
\label{eqd7}
\nu^N(t,dx)=\left(\frac{N-1}{N}\right)^{A^N_t}\mu^N(t,dx),
\end{eqnarray*}
where $A^N(t)=\#\bigcup_{i=1}^N\{\tau^i_n,\ \ 0\leq\tau^i_n\leq t\}$ denotes the total number of jumps before time $t$.

\begin{lemme}
 \label{led1}
 The sequence of measure processes $\left(\nu^N(.,dx)\right)_N$ converges in law to $P_{\mu_0}(X. \in dx)$ in the Skorokhod topology $D([0,T],{\cal M}(]0,1[))$.
\end{lemme}
\textit{Proof of Lemma \ref{led1} :}
Applying the It\^o's formula to the semimartingale $\nu^N(t,\psi)$, we deduce from \eqref{eqd5} that
\begin{align*}
 \nu^N(t,\psi)=\nu^N(0,\psi)+\int_0^t{\nu^N(s\minus,\psi'q+\psi'')ds}&+\int_0^t{\left(\frac{N-1}{N}\right)^{A^N_{s\minus}}d{\cal M}^c(\psi,s)}\\&+\sum_{0\leq\tau_n\leq t}(\nu^N(\tau_n,\psi)-\nu^N(\tau_n\minus,\psi)),
\end{align*}
Where we have
\begin{multline*}
 \nu^N(\tau_n,\psi)-\nu^N(\tau_n\minus,\psi)=\left(\frac{N-1}{N}\right)^{A^N_{\tau_n}}\left(\mu^N(\tau_n,\psi)-\mu^N(\tau_n\minus,\psi)\right)\\
+\mu^N(\tau_n\minus,\psi)\left(\left(\frac{N-1}{N}\right)^{A^N_{\tau_n}}-\left(\frac{N-1}{N}\right)^{A^N_{\tau_n\minus}}\right),
\end{multline*}
with
\begin{eqnarray*}
 \label{eqd12}
 \mu^N(\tau_n,\psi)-\mu^N(\tau_n\minus,\psi)=\frac{1}{N-1}\mu^N(\tau_n\minus,\psi)+{\cal M}^j(\psi,\tau_n)-{\cal M}^j(\psi,\tau_n\minus)
\end{eqnarray*}
and
\begin{eqnarray*}
 \label{eqd10}
 \left(\frac{N-1}{N}\right)^{A^N_{\tau_n}}-\left(\frac{N-1}{N}\right)^{A^N_{\tau_n\minus}}=-\frac{1}{N-1}\left(\frac{N-1}{N}\right)^{A^N_{\tau_n}},
\end{eqnarray*}
then
\begin{eqnarray*}
 \nu^N(\tau_n,\psi)-\nu^N(\tau_n\minus,\psi)&=&\left(\frac{N-1}{N}\right)^{A^N_{\tau_n}}\left({\cal M}^j(\psi,\tau_n)-{\cal M}^j(\psi,\tau_n\minus) \right).\\
    &=&\frac{N-1}{N}\left(\frac{N-1}{N}\right)^{A^N_{\tau_n\minus}}\left({\cal M}^j(\psi,\tau_n)-{\cal M}^j(\psi,\tau_n\minus) \right)
\end{eqnarray*}
That implies
\begin{multline*}
 \nu^N(t,\psi)-\nu^N(0,\psi)-\int_0^t{\nu^N(s\minus,\psi'q+\frac{1}{2}\psi'')}ds=\int_0^t{\left(\frac{N-1}{N}\right)^{A^N_{s\minus}}d{\cal M}^c(\psi,s)}\\
+\frac{N-1}{N}\sum_{0\leq\tau_n\leq t}{\left(\frac{N-1}{N}\right)^{A^N_{\tau_n\minus}}\left({\cal M}^j(\psi,\tau_n)-{\cal M}^j(\psi,\tau_n\minus)\right)}
\end{multline*}
We deduce that for all smooth functions $\Psi(t,x)$ vanishing at the boundary
\begin{multline}
 \nu^N(t,\Psi(t,.))-\nu^N(0,\Psi(0,.))-\int_0^t{\nu^N(s\minus, \frac{\partial \Psi(s,.)}{\partial s}+\frac{\partial\Psi(s,.)}{\partial x}q+\frac{1}{2}\frac{\partial^2\Psi(s,.)}{\partial x^2})}ds\\
=\int_0^t{\left(\frac{N-1}{N}\right)^{A^N_{s\minus}}d{\cal M}^c(\Psi(s,.),s)}\\
+\frac{N-1}{N}\sum_{0\leq\tau_n\leq t}{\left(\frac{N-1}{N}\right)^{A^N_{\tau_n\minus}}\left({\cal M}^j(\Psi(\tau_n,.),\tau_n)-{\cal M}^j(\Psi(\tau_n\minus,.),\tau_n\minus)\right)}
\label{eqd14}
\end{multline}
Because $\left(\frac{N-1}{N}\right)^{A^N_{s\minus}}\leq 1$ a.s. and by the Doob's inequality , we have
\begin{eqnarray}
\label{eqd18}
 E\left(\sup_{t\in[0,T]} \left| \int_0^t{\left(\frac{N-1}{N}\right)^{A^N_{s\minus}}d{\cal M}^c(\Psi(t,.),s)} \right|^2\right)
\leq \frac{1}{N} T \|\frac{\partial \Psi}{\partial x}\|_{\infty}^2.
\end{eqnarray}
Note that the jumps of the martingale ${\cal M}^j$ are smaller than $\frac{2}{N}\|\Psi\|_{\infty}$, then
\begin{align*}
 E\Bigg[\sum_{0\leq \tau_n \leq T}{\left(\frac{N-1}{N}\right)^{2 A_{\tau_n\minus}}}&{\left({\cal M}^j(\Psi(\tau_n,.),\tau_n)-{\cal M}^j(\Psi(\tau_n\minus,.),\tau_n\minus)\right)^2}\Bigg]\\
&\leq \frac{4}{N^2}\|\Psi\|_{\infty}^2 E\left[\sum_{0\leq \tau_n \leq T}{\left(\frac{N-1}{N}\right)^{2 A_{\tau_n\minus}}}\right]\\
&\leq \frac{4}{N}\|\Psi\|_{\infty}^2
\end{align*}
By the Doob's inequality, we have then
\begin{equation}
\label{eqd19}
E\left(\sup_{t\in[0,T]} \left|\sum_{0\leq\tau_n\leq t}{\left(\frac{N-1}{N}\right)^{A^N_{\tau_n\minus}}\left({\cal M}^j(\Psi(\tau_n,.),\tau_n)-{\cal M}^j(\Psi(\tau_n\minus,.),\tau_n\minus)\right)} \right|^2\right)
\leq \frac{4}{N}\|\Psi\|_{\infty}^2,
\end{equation}

Define $\Psi(s,x)=P_{t-s}f(x)$, where $f\in C^{\infty}([0,1])$ vanishes on $\{0,1\}$, and $(P_t)$ is the semigroup associated with the diffusion $X$ defined by \eqref{eq1}. From Kolmogorov's equation (see \cite[Proposition 1.5 p.9]{et1}),
\begin{eqnarray*}
	\label{eq24}
	\frac{\partial}{\partial s}\Psi(s,x)+\frac{1}{2}\Delta\Psi(s,x)+q(x)\nabla\Psi(s,x)=0.
\end{eqnarray*}
We deduce from \eqref{eqd14}, \eqref{eqd18} and \eqref{eqd19}, that
\begin{eqnarray*}
E\left( \sup_{t\in[0,T]} \left| \nu^N(t,f)-\int_0^1{P_{t}f(x)}d\mu^N(0,x) \right|^2\right) \leq \frac{1}{N} C(f),
\end{eqnarray*}
where $C(f)$ is a positive constant, which only depends on $f$.
For each map $g\in C^{\infty}([0,1])$, one can set $f_r(x)=\gamma_r(x)g(x)$, with $r>0$, where $\gamma_r\in C^\infty([0,1])$ is equal to $1$ on $]2r, 1-2r[$ and vanishes on $]0,r[\cup]1-r,1[$. Then
\begin{eqnarray}
 \left| \nu^N(t,g)-\int_0^1{P_{t}g(x)}d\mu^N(0,x) \right| &\leq& \left| \nu^N(t,f_r)-\int_0^1{P_{t}f_r(x)}d\mu^N(0,x) \right| \label{eqd15}\nonumber\\
&&+ \left| \nu^N(t,(1-\gamma_r)g) \right|\label{eqd16}\\
&&+\left| \int_0^1{P_{t}\left((1-\gamma_r)g\right)(x)}d\mu^N(0,x) \right|\label{eqd17},
\end{eqnarray}
where \eqref{eqd16} (see Proposition \ref{pr4}) and \eqref{eqd17} are going to $0$ when $r$ tends to $0$, uniformly in $N$, and $f_r\in C^{\infty}([0,1])$ vanishes on $\{0,1\}$. Then
\begin{eqnarray*}
 \sup_{t\in[0,T]} \left| \nu^N(t,g)-\int_0^1{P_{t}g(x)}d\mu^N(0,x) \right| \FlechBH{\longrightarrow}{N\rightarrow\infty}{L^2} 0.
\end{eqnarray*}
In particular, $\nu^N(.,dx)$ converges in law to $P_{\mu_0}(X. \in dx)$. $\Box$

\paragraph{}
Let us conclude the proof of Proposition \ref{pr5}.
From Lemma \ref{led1},
\begin{equation*}
(\nu^N(t,]0,1[),\nu^N(t,dx))_{t\in[0,T]} \FlechBH{\longrightarrow}{N\rightarrow\infty}{law}(\mathbb{P}_{\mu(0,dx)}(X_t\in]0,1[),\mathbb{P}_{\mu(0,dx)}(X_t\in dx))_{t\in[0,T]}
\end{equation*}
in the Skorokhod topology $D([0,T],\mathbb{R}\times{\cal M}_1(]0,1[))$. That means
\begin{equation*}
 \left(\left(\frac{N-1}{N}\right)^{A^N_{t}},\nu^N(t,dx)\right)_{t\in[0,T]}\FlechBH{\longrightarrow}{N\rightarrow\infty}{law}(\mathbb{P}_{\mu(0,dx)}(X_t\in]0,1[),\mathbb{P}_{\mu(0,dx)}(X_t\in dx))_{t\in[0,T]}.
\end{equation*}
		The process $\mathbb{P}_{\mu(0,dx)}(X.\in]0,1[)$ never vanishing almost surely and the limit process being continuous almost surely, we have
\begin{equation*}
 \left(\left(\frac{N-1}{N}\right)^{A^N_{t}}\nu^N(t,dx)\right)_{t\in[0,T]}\FlechBH{\longrightarrow}{N\rightarrow\infty}{law}\left(\mathbb{P}_{\mu(0,dx)}(X_t\in dx)/\mathbb{P}_{\mu(0,dx)}(X_t\in]0,1[)\right)_{t\in[0,T]}
\end{equation*}
 in the Skorokhod topology $D([0,T],{\cal M}_1(]0,1[))$. That means
\begin{equation*}
(\mu^N(t,dx))_{t\in[0,T]} \FlechBH{\longrightarrow}{N\rightarrow\infty}{law} (\mathbb{P}_{\mu(0,dx)}(X_t\in dx|X_t\neq\partial))_{t\in[0,T]}
\end{equation*}
in the Skorokhod topology $D([0,T],{\cal M}_1(]0,1[))$.
		The proof of Proposition \ref{pr5} is then complete.

\subsection{Existence and convergence of the empirical stationary measures ${\cal X}^N$}
\label{parts23}
For each $N\geq2$, we say that the interacting particle process $(X^1,...,X^N)$ is exponentially ergodic, if there exists a probability measure $M^N$ on $]0,1[^N$ such that,
\begin{equation}
 ||P_x((X^1_t,...,X^N_t)\in .)-M^N||_{TV}\leq C(x)\rho^t,\ \forall x\in]0,1[^N,\ \forall t\in\mathbb{R}_+,
\end{equation}
where $C(x)$ is finite, $\rho<1$ and $||.||_{TV}$ is the total variation norm. In particular, $M^N$ is a stationary measure for the process $(X^1,...,X^N)$. When $M^N$ exists, we denote by ${\cal X}^N$ the empirical stationary measure associated with $M^N$, that is a random probability which is distributed as $\frac{1}{N}\sum_{i=1}^N{\delta_{x_i}}$, where $(x^1,...,x^N)$ is a random vector in $]0,1[^N$ distributed with respect to $M^N$.

In a first time, we prove that for all $N\geq 2$, the interacting particle process with $N$ particles $(X^1,...,X^N)$ associated with $X$ is exponentially ergodic. We conclude by proving that $({\cal X}^N)_N$ converges in law to the unique QSD of $X$.

	\subsubsection{Exponential ergodicity}
Here $N\geq 2$ is fixed. We are interested in proving the following result, which is the first part of Theorem \ref{th2}

\begin{proposition}
\label{pr9}
 The interacting particle process $(X^1,...,X^N)$ with law $\mathbb{P}^{ipp}$ is exponentially ergodic.
\end{proposition}
\textit{Proof of Proposition \ref{pr9} : }We focus on the $1$-skeleton of the interacting particle process with $N$ particles, which is the Markov chain $(X^1_n,...,X^N_n)_{n\in\mathbb{N}}$.  Thanks to \cite[Theorem 5.3 p.1681]{do1}, exponential ergodicity of $(X^1,...,X^N)$ will be obtained as soon as the associated $1$-skeleton is geometrically ergodic, which means that it exists a probability measure $\pi^N$ on $]0,1[^N$ such that
\begin{eqnarray*}
 ||P_x((X^1_n,...,X^N_n)\in .)-\pi^N||_{TV}\leq C_0(x)\rho_0^n,\ \forall x\in]0,1[^N,\ \forall n\in\mathbb{N},
\end{eqnarray*}
where $C_0(x)$ is finite and $\rho_0<1$.

To prove the geometrical ergodicity of the 1-skeleton, let us introduce the following definition:
\begin{definition}
 $C\subset]0,1[^N$ is said to be a \textbf{small set} for the Markov chain $(X^1_n,...,X^N_n)_{n\in\mathbb{N}}$ if, for some nontrivial probability measure $\vartheta$ and some $n\geq 1$, $\epsilon>0$, the $n$-step transition probability kernel $P^n(x,A):=\mathbb{P}^{ipp}_x((X^1_n,...,X^N_n)\in A)$ satisfies, for all $x\in C$,
		\begin{equation*}
		 P^n(x,A)\geq\epsilon\vartheta(A),\ A\in{\cal B}(]0,1[^N).
		\end{equation*}
\end{definition}

\begin{lemme}
\label{le4}
All compact set $C=[r,1-r]^N$ (with $r>0$) is a small set for the $1$-skeleton. Moreover, $\exists \kappa>0$ so that
\begin{eqnarray}
	\label{eq27}
	\sup_{x\in C}{E_x(\kappa^{\tau'_C})}<\infty,
\end{eqnarray}
where $\tau'_C$ is the return time to $C$.
\end{lemme}
\textit{Proof of Lemma \ref{le4}:} Fix $r>0$ and let ${\cal F}$ be the event ``the process $(X^1,...,X^N)$ has no jumps between times $0$ and $1$''.  Define $p=\inf_{x\in[r,1-r]^N}\mathbb{P}^{ipp}_x({\cal F})$. Thanks to the coupling with $(Y^1,...,Y^N)$, we have $p>0$. Conditionally to the event $\cal F$, the particles of the interacting particle process are independent from each other.

Let us study $\vartheta^1(dx)=\mathbb{P}^{ipp}_{(x_1,...,x_N)}(X^1_1\in dx\ \mbox{and}\ {\cal F})$. The law of $X^1$ conditionally to $\cal F$ is the same as the law of $X$ conditioned to not jump, because, given this last event, the process $X^1$ doesn't depend on the other particles. Thus the probability of ``$X^1_1\in dx\ \mbox{and}\ {\cal F}$'' is $\mathbb{P}_{x_1}(X_1\in dx)$. The law of $X_1$ has a density $p_1(x_1,y)$ with respect to the Lebesgue's measure and $p_1(x_1,y)$ depends continuously on $x_1$ and $y$. It only vanishes when $y=0$ or $1$. Then
\begin{equation*}
 \inf_{(x_1,y)\in[r,1-r]\times[r,1-r]}{p_1(x_1,y)}>0
\end{equation*}
Denoting this minimum by $\epsilon'$, we have, for all $x_1\in[r,1-r]$, $\vartheta^1(dx)\geq \epsilon' \mathbf{1}_{[r,1-r]}(x) dx$.

Conditionally to $\cal F$, the particles are independent from each other, so that
\begin{eqnarray*}
\mathbb{P}^{ipp}_{(x_1,...,x_N)}((X^1_1,...,X^N_1)\in dy_1...dy_N|{\cal F})=\prod_{i=1}^N{\mathbb{P}^{ipp}_{(x_1,...,x_N)}(X^i_1\in dy_i|{\cal F})},
\end{eqnarray*}
where $\mathbb{P}^{ipp}_{(x_1,...,x_N)}(X^i_1\in dy_i|{\cal F})$ is greater than $\mathbb{P}^{ipp}_{(x_1,...,x_N)}(X^i_1\in dy_i\ \mbox{and}\ {\cal F})$ and then greater than $\vartheta^1(dy_i)$. Finally, we have
\begin{eqnarray*}
 \mathbb{P}^{ipp}_{(x_1,...,x_N)}((X^1_1,...,X^N_1)\in dy_1...dy_N|{\cal F})\geq {\epsilon'}^N \mathbf{1}_{[r,1-r]^N}(y_1,...,y_N) dy_1...dy_N,
\end{eqnarray*}
so that $[r,1-r]^N$ is a small set.

For all $x\in]0,1[^N$, and all $n\geq1$, the probability of being in $C$ at time $n+1$ starting from $x$ at time $n$
is bounded below by the probability $p_C>0$ for $(Y^1,...,Y^N)$ to enter $C$ at time $n+1$, starting from $0$ at time $n$. Hence, at each time $n\geq1$, $(X^1,...,X^N)$ returns to $C$ at time $n+1$ with a probability greater than $p_C>0$.
That implies that the return time to $C$ for the 1-skeleton of the interacting particle process with $N$ particles is bounded above by a time of geometrical law, independent of the starting point $x\in]0,1[^N$, and then satisfies condition \eqref{eq27}.
$\Box$

\paragraph{}
The chain $(X^1_n,...,X^N_n)_{n\in\mathbb{N}}$ is aperiodic. Moreover, if the Lebesgue measure of a subset $A\subset]0,1[^N$ is strictly positive, then $\mathbb{P}^{ipp}_x(\tau_A<\infty)>0$ for all $x\in]0,1[^N$, where $\tau_A$ is the first hitting time on $A$ for the chain $(X^1_n,...,X^N_n)_{n\in\mathbb{N}}$. Thanks to \cite[Theorem 2.1 p.1673]{do1}, if such a Markov chain has a small set which satisfies \eqref{eq27}, then it is geometrically ergodic. As a consequence, Lemma \ref{le4} allows us to conclude the proof of Proposition \ref{pr9}.

	\subsubsection{Convergence to the QSD}
 We are interested in proving the following result, which is the second part of Theorem \ref{th2}
\begin{proposition}
\label{pr10}
 The sequence of random measures $\left({\cal X}^N\right)_{N\geq2}$ converges in law to the deterministic measure $\nu$, QSD of the process $X$.
\end{proposition}
\textit{Proof of Proposition \ref{pr10} : }
 For each $r\in]0,1/4[$, we define $\gamma_r$ as a non-negative continuous bounded function from $]0,1[$ to $\mathbb{R}$, equal to $1$ on $D_{2r}$ and equal to $0$ on $D_{r}^c$. We have 
\begin{equation*}
{\cal X}^N(D_r^c)\leq {\cal X}^N(1-\gamma_r),\ \forall r\in]0,1/4[. 
\end{equation*}
	Thanks to Proposition \ref{pr9}, the sequence of random measures $(\mu^N(t,dx))_{N\geq2}$ converges in law to ${\cal X}^N$ when $t$ tends to $+\infty$. That implies
\begin{equation*}
 E\left(\mu^N(t,1-\gamma_r)\right)\FlechBH{\longrightarrow}{t\rightarrow+\infty}{} E\left({\cal X}^N(1-\gamma_r)\right),\ \forall r\in]0,1/4[.
\end{equation*}
We denote by $\mu'^N(t,dx)$ the empirical measure of $(Y^1_t,...,Y^N_t)$. Let us choose $\gamma_r$ monotone on $]r,2r[$ and $]1-2r,1-r[$. From the coupling inequality, 
\begin{eqnarray*}
 \mu^N(t,1-\gamma_r)\leq\mu'^N(t,1-\gamma_r)
\end{eqnarray*}
for all $t\in\mathbb{R}_+$ and $r>0$. Then
\begin{eqnarray*}
 E\left(\mu^N(t,1-\gamma_r)\right)\leq E\left(\mu'^N(t,1-\gamma_r)\right),
\end{eqnarray*}
which tends to $0$ when 
$t$ tends to $+\infty$ and $r$ to $0$, uniformly in $N$. As a consequence,
\begin{equation*}
 E({\cal X}^N)(1-\gamma_r) \FlechBH{\longrightarrow}{r\rightarrow0}{} 0,
\end{equation*}
where $E({\cal X}^N)$ is the deterministic measure defined by $E({\cal X}^N)(A)=E({\cal X}^N(A))$, for all measurable set $A$.
That yields
\begin{equation*}
 E({\cal X}^N)(D_r^c) \FlechBH{\longrightarrow}{r\rightarrow0}{} 0,
\end{equation*}
uniformly in $N$.
The family of intensity measures $(E({\cal X}^N))_{N}$ is then tight. This is a sufficient condition for 
	the family of random variables $({\cal X}^N)$ to be tight, as shown in \cite[Corollary 2.2]{ja2}. We conclude that it exists a sub-sequence $({\cal X}^{\phi(N)})$ which converges in law to a random probability measure $\cal X$.

	Choose $\mu^N(0,dx)={\cal X}^N(dx)$. The non-degeneracy property is fulfilled. Thanks to Proposition \ref{pr5},
\begin{equation*}
 (\mu^{\phi(N)}(t,dx))_{t\in[0,T]} \FlechBH{\longrightarrow}{N\rightarrow\infty}{law} (\mathbb{P}_{\cal X}(X_t\in dx|X_t\neq \partial))_{t\in[0,T]},\ \forall T>0
\end{equation*}
in the Skorokhod topology $D([0,T],{\cal M}_1(]0,1[))$. The limiting process $(\mathbb{P}_{\cal X}(X_t\in dx|X_t\neq \partial))_{t\in[0,T]}$ being almost surely continuous,
\begin{equation*}
 \mu^{\phi(N)}(t,dx) \FlechBH{\longrightarrow}{N\rightarrow\infty}{law} \mathbb{P}_{\cal X}(X_t\in dx|X_t\neq \partial),\ \forall t>0,
\end{equation*}
with respect to the weak topology of ${\cal M}_1(]0,1[)$.
	By stationarity, the random probability measures $\mu^{\phi(N)}(t,dx)$ and ${\cal X}^{\phi(N)}$ have the same law. Making $N$ tend to $\infty$, we deduce that $\mathbb{P}_{\cal X}(X_t\in dx|X_t\neq \partial)$ and ${\cal X}$ have the 
	same law too. This looks like a QSD, but ${\cal X}$ is \textit{a priori} a random measure and we need the following result to conclude.
	\begin{lemme}
		\label{le3}
		For all $m\in{\cal M}_1(]0,1[)$,
		\begin{eqnarray*}
		\label{eq39}
		\lim_{t\rightarrow\infty}{P_{m}(X_t\in A | \tau>t)}=\nu(A),
		\end{eqnarray*}
	where $\tau$ is the killing time of the process $X$ and $\nu$ its unique QSD.
	\end{lemme}
	\textit{Proof of Lemma \ref{le3} :} Let $m$ be a probability measure on $]0,1[$. $\exists \lambda_0>0$, $\phi_0$ and $\widetilde{\phi}_0$ two continuous maps vanishing on $0$ and $1$ such that, for all $x\in ]0,1[$ (see R.G. Pinsky's explanations \cite[Hypotheses 2 and 3]{pi1}):
	\begin{align}
		\label{eq40}
		\lim_{t\rightarrow\infty}{e^{\lambda_0 t}P_x(\tau>t)}=C_1 \phi_0(x)&,\\
		\label{eq41}
		\lim_{t\rightarrow\infty}{e^{\lambda_0 t}P_x(X_t\in A,\ \tau>t)}=&\;C_2\phi_0(x)\int_A{\widetilde{\phi}_0(y)}dy.
	\end{align}
	Here $e^{\lambda_0 t}P_x(\tau>t)$ is uniformly bounded above in the variables $t$ and $x$ (see \cite[Proof of the equality 7.2, p27]{ca1}),
	then, by dominated convergence, one can integrate with respect to $m$ under the limit in \eqref{eq40},
	\begin{eqnarray*}
		\label{eq42}
		\lim_{t\rightarrow\infty}{e^{\lambda_0 t}P_m(\tau>t)}=C_1 \int_D{\phi_0(x)}m(dx).
	\end{eqnarray*}
	The same holds for \eqref{eq41}:
	\begin{eqnarray*}
		\label{eq43}
		\lim_{t\rightarrow\infty}{e^{\lambda_0 t}P_m(X_t\in A,\ \tau>t)}=C_2\int_D{\phi_0(x)\int_A{\widetilde{\phi}_0(y)}dy}\ m(dx).
	\end{eqnarray*}
		Then, by Fubini's Theorem,
		\begin{eqnarray*}
			\label{eq44}
			\lim_{t\rightarrow\infty}{\frac{P_m(X_t\in A,\ \tau>t)}{P_m(\tau>t)}}=\frac{C_2}{C_1}\int_A{\widetilde{\phi}_0(y)}dy,
		\end{eqnarray*}
		that is, from (\ref{eq40}) and (\ref{eq41}) with $A=]0,1[$,
		\begin{eqnarray*}
			\label{eq45}
			\lim_{t\rightarrow\infty}{\frac{P_m(X_t\in A,
							\ \tau>t)}{P_m(\tau>t)}}=\frac{\int_A{\widetilde{\phi}_0(y)}dy}{\int_D{\widetilde{\phi}_0(y)}dy}
		\end{eqnarray*}
		which is nothing else but $\nu(A)$ (see \cite[Proposition 1.10]{pi1}). $\Box$

\paragraph{}
	Thanks to Lemma \ref{le3}, $\mathbb{P}_{\cal X}(X_t\in dx|X_t\neq \partial)$ converges almost surely to the Yaglom limit when $t\rightarrow+\infty$, and so do $\cal X$.
	Finally, $\cal X$ is the unique QSD of the process and the proof is complete. $\Box$

\section{Numerical applications}

\subsection{The logistic case}
\label{par5}
We apply our result to the logistic Feller diffusion with values in $]0,+\infty[$, driven by the stochastic differential equation
\begin{equation}
 \label{eq77}
 dZ_t=\sqrt{Z_t}dB_t+(rZ_t-cZ_t^2)dt,\ Z_0=z>0,
\end{equation}
and killed when it hits $0$. Here $B$ is a $1$-dimensional Brownian motion and $r,c$ are two positive constants.

We define $\mathbb{P}^0$ as the law of $2\sqrt{Z.}$, which is killed at $0$ and satisfies the SDE
\begin{equation*}
 dX_t=dB_t-\left(\frac{1}{2X_t}-\frac{rX_t}{2}+\frac{cX_t^3}{4}\right)dt,\ X_0=x\in]0,+\infty[.
\end{equation*}
For each $\epsilon>0$, we define the law $\mathbb{P}^{\epsilon}$ and denote its QSD by $\nu_{\epsilon}$.

As proved in \cite{ca1}, (H1) and (H3) are fulfilled in this case. Thanks to Theorem \ref{th1} and denoting by $\nu$ the Yaglom limit associated with $\mathbb{P}^0$, we have
\begin{equation*}
\nu_{\epsilon}\FlechBH{\longrightarrow}{\epsilon\rightarrow0}{law}\nu.
\end{equation*}
In the numerical simulations below, we set $\epsilon$ equal to $0.001$.

\begin{figure}[htbp]
\begin{center}
 \input{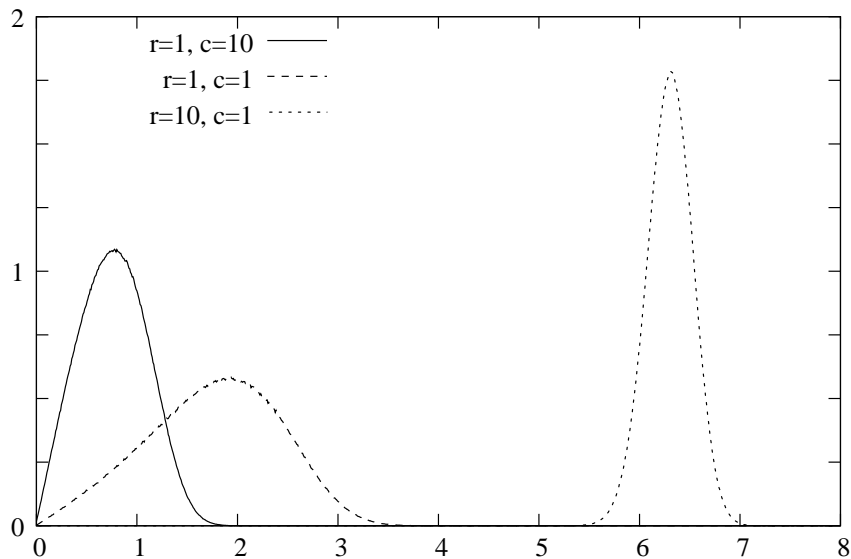}
\caption{$E({\cal X}^N)$ for the diffusion \eqref{eq77}, with different values of $r$ and $c$}
\label{fig10}
\end{center}
\end{figure}

By Theorem \ref{th2}, we have
\begin{equation*}
 {\cal X}^N\FlechBH{\longrightarrow}{N\rightarrow+\infty}{} \nu_{\epsilon},
\end{equation*}
where ${\cal X}^N$ is the empirical measure of the system studied in Section \ref{par2}.
In the numerical simulations, we set $N=1000$ and, because of the randomness of ${\cal X}^N$, we approximate $E({\cal X}^N)$ using the Ergodic theorem: we compute $\frac{1}{10000}\sum_{t=1}^{10000}{\mu_N(t,dx)}$.
The graphic below (see Figure \ref{fig10}) shows this approximation for different values of $r$ and $c$.

As it could be wanted for, greater is $c$, closer is the support of the QSD to $0$. We thus numerically describe the impact of the linear and quadratic terms on the QSD.




\subsection{The Wright-Fisher case}

We illustrate the result of Section \ref{par3} by an application to the Wright-Fisher diffusion with values in $]0,1[$ conditioned to be killed at $0$. This diffusion is driven by the SDE 
\begin{equation*}
 dZ_t=\sqrt{Z_t(1-Z_t)}dB_t-Z_t dt,\ Z_0=z\in]0,1[,
\end{equation*}
and killed when it hits $0$ ($1$ is never reached). In \cite{hu1}, the author proves that the QSD of this process exists and has the density $2-2x$ with respect to the Lebesgue measure.

Define $\mathbb{P}^0$ as the law of $X_.=\arccos(1-2 Z_.)$, where $Z$ is defined as above. $\mathbb{P}^0$ is the law of the diffusion with values in $]0,\pi[$, driven by the SDE
\begin{equation*}
 dX_t=dB_t-\frac{1-2\cos X_t}{2\sin X_t } dt,\ X_0=x\in]0,\pi[,
\end{equation*}
killed when it hits $0$ ($\pi$ is never reached). For all $\epsilon\in]0,\pi/2[$, define $\mathbb{P}^{\epsilon}$ and $\nu_{\epsilon}$ as in Section \ref{par3}.

The drift of the diffusion is $q(x)=\frac{1-2\cos X_t }{2\sin X_t}$, $\forall x\in]0,\pi[$. Let us show that it satisfies the hypotheses of Theorem \ref{th10}.

We have, $\forall x\in]0,\pi[$,
\begin{eqnarray*}
 q(x)^2-q'(x)&=&\frac{8 \cos^2 x +2 \cos x +2 \cos x  +1}{4\sin^2 x }+1
\end{eqnarray*}
which is positive and tends to $+\infty$ both in $0+$ and $\pi-$.
It implies that hypothesis (HH1) is fulfilled.

For all $x\in]0,\pi[$,
\begin{eqnarray*}
 Q(x)&=&\int_{\pi/2}^{x}{q(y)dy}\\
     &=&\frac{1}{2}\ln\left|\tan\frac{x}{2}\right|+\ln\left|\sin x\right|\\
     &=&\ln\left(2\left(\sin\frac{x}{2}\right)^{\frac{3}{2}}\left(\cos\frac{x}{2}\right)^{\frac{1}{2}}\right).
\end{eqnarray*}
That implies
\begin{equation*}
e^{-Q(x)}=\frac{1}{2\left(\sin\frac{x}{2}\right)^{\frac{3}{2}}\left(\cos\frac{x}{2}\right)^{\frac{1}{2}}}
\end{equation*}
then
\begin{equation*}
  x e^{-Q(x)} \;\FlechBH{\sim}{0+}{}\; \sqrt{\frac{2}{x}}
\end{equation*}
and
\begin{equation*}
 (\pi-x)  e^{-Q(x)} \;\FlechBH{\sim}{\pi-}{}\; \sqrt{\frac{\pi-x}{2}}.
\end{equation*}
Finally, hypotheses (HH2) and (HH3) are satisfied and Theorem \ref{th10} can be applied.

\begin{figure}[htbp]
\begin{center}
\input{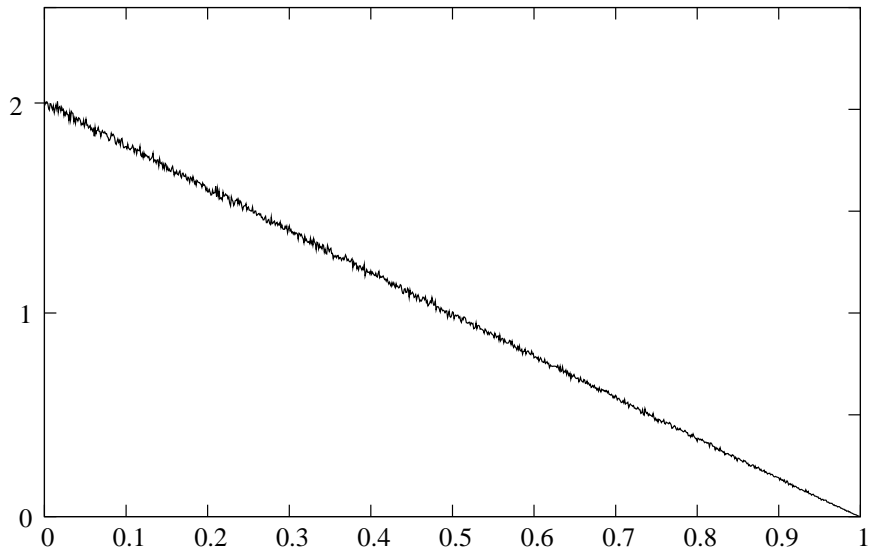}
\caption{Wright Fisher case, $\epsilon=0.001$ and $N=1000$}
\label{fig4}
\end{center}
\end{figure}

In the following numerical simulation (see Figure \ref{fig4}), we set $\epsilon=0.001$ and $N=1000$. We compute $E({\cal X}^N)$, which is an approximation of $\nu_{\epsilon}$, and then of $\nu$, with the method used in the logistic case (see Part \ref{par5}).

The simulation is very close to the QSD $(2-2x)dx$, which shows the efficiency of the method.

\begin{remarque}\upshape
 In the simulations, we have chosen to simulate a system with $N=1000$ particles, because of empirical consideration. The question of convergence speed will be studied in further works.
\end{remarque}

\paragraph{Acknowledgments} I am extremely grateful to my Ph.D. supervisor Sylvie M\'el\'eard for his carefully and essential help on the form and the content of my first paper. I would like to thank Pierre Collet for his explanations on some spectral theory tools and Krzysztof Burdzy for his indication on the proof of the interacting particle process existence.

\end{document}